\documentclass[11pt,reqno]{amsart}
\usepackage{amsmath,epsfig,graphicx,color}
\usepackage{ifthen}
\usepackage{amssymb,latexsym}
\usepackage{subfigure}
\usepackage{hyperref}
\hypersetup{
urlcolor=black,
  menucolor=black,
  citecolor=black,
  anchorcolor=black,
  filecolor=black,
  linkcolor=black,
  colorlinks=true,
}
\newcommand{\rep}{representation}
\newcommand{\ld}{large deviation}
\renewcommand{\H}{{\mathbf H}}

\newcommand{\cmt}{continuous mapping theorem}
\newcommand{\fct}{function}
\newcommand{\slvary}{slowly varying}
\newcommand{\regvar}{regular variation}
\newcommand{\regvary}{regularly varying}
\newcommand{\st}{such that}
\newcommand{\stas}{\stackrel{\rm a.s.}{\rightarrow}}
\newcommand{\std}{\stackrel{\rm d}{\rightarrow}}
\newcommand{\stp}{\stackrel{\P}{\rightarrow}}
\newcommand{\la}{\lambda}
\newcommand{\ds}{distribution}
\newcommand{\beao}{\begin{eqnarray*}}
\newcommand{\eeao}{\end{eqnarray*}}
\newcommand{\beam}{\begin{eqnarray}}
\newcommand{\eeam}{\end{eqnarray}}

\definecolor{darkblue}{rgb}{.1, 0.1,.8}
\definecolor{darkgreen}{rgb}{0,0.8,0.2}
\definecolor{darkred}{rgb}{.8, .1,.1}
\newcommand{\wt}{\widetilde}
\newcommand{\bco}{\begin{corrolary}}
\newcommand{\eco}{\end{corrolary}}

\newcommand{\E}{\mathbb{E}}
\renewcommand{\P}{\mathbb{P}}
\newcommand{\1}{\mathbf{1}}
\newcommand{\R}{\mathbb{R}}
\newcommand{\N}{\mathbb{N}}

\newcommand{\Z}{\mathbb{Z}}

\newcommand{\Frechet}{Fr\'{e}chet }

\DeclareMathOperator{\e}{e}

\newcommand{\X}{{\mathbf X}}

\newcommand{\M}{{\mathbf M}}

\newcommand{\dint}{\,\mathrm{d}}

\newcommand{\twonorm}[1]{\|#1\|_2}
\newcommand{\inftynorm}[1]{\|#1\|_\infty}
\newcommand{\frobnorm}[1]{\|#1\|_F}
\newcommand{\vep}{\varepsilon}
\newcommand{\nto}{n \to \infty}
\newcommand{\xto}{x \to \infty}

\newcommand{\rhs}{right-hand side}
\newcommand{\ts}{time series}
\newcommand{\tsa}{\ts\ analysis}
\newcommand{\fidi}{finite-dimensional distribution}
\newcommand{\rv}{random variable}

\newcommand{\diag}{\operatorname{diag}}
\newcommand{\slln}{strong law of large numbers}
\newcommand{\clt}{central limit theorem}

\newcommand{\ex}{{\rm e}\,}

\def\tag{\refstepcounter{equation}\leqno }

\newtheorem{lemma}{Lemma}[section]

\newtheorem{theorem}[lemma]{Theorem}
\newtheorem{proposition}[lemma]{Proposition}

\newtheorem{corollary}[lemma]{Corollary}
\newtheorem{example}[lemma]{Example}

\newtheorem{remark}[lemma]{Remark}

\newcommand{\cid}{\stackrel{d}{\rightarrow}}
\newcommand{\cip}{\stackrel{\P}{\rightarrow}}
\newcommand{\civ}{\stackrel{v}{\rightarrow}}

\newcommand{\A}{{\mathbf A}}

\newcommand{\as}{{\rm a.s.}}

\newcommand{\evt}{extreme value theory}

\newcommand{\pp}{point process}
\newcommand{\con}{convergence}
\newcommand{\seq}{sequence}
\newcommand{\ms}{measure}
\newcommand{\asy}{asymptotic}
\begin{document}
\today
\title[Extreme value analysis for the sample autocovariance matrices of time series]
{Extreme value analysis for the sample autocovariance matrices of heavy-tailed multivariate time series}
\thanks{Richard Davis
was supported by ARO MURI grant W911NF-12-1-0385. Thomas Mikosch's and Johannes Heiny's research is partly supported by the Danish Research Council Grant DFF-4002-00435 ``Large random matrices with heavy tails and dependence''.}
\author[Richard A. Davis]{Richard A. Davis}
\author[Johannes Heiny]{Johannes Heiny}
\author[Thomas Mikosch]{Thomas Mikosch}
\author[Xiaolei Xie]{Xiaolei Xie}
\address{Department of Statistics,
Columbia University,
1255 Amsterdam Ave.
New York, NY 10027, U.S.A.}
\email{rdavis@stat.columbia.edu\,, www.stat.columbia.edu/$\sim$rdavis}
\address{Department  of Mathematics,
University of Copenhagen,
Universitetsparken 5,
DK-2100 Copenhagen,
Denmark}
\email{johannes.heiny@math.ku.dk}
\email{mikosch@math.ku.dk\,, www.math.ku.dk/$\sim$mikosch}
\email{xie@math.ku.dk}
\begin{abstract}
We  provide some  asymptotic theory for the largest eigenvalues of a sample covariance matrix
of a $p$-dimensional \ts\ where the dimension $p=p_n$ converges to infinity when the sample size $n$ increases.
We give a short overview of the literature on the topic both in the light- and heavy-tailed cases when the data have
finite (infinite) fourth moment, respectively.
Our main focus is on the heavy-tailed case. In this case, one has a theory for the \pp\ of the normalized eigenvalues
of the sample covariance matrix in the iid case but also when rows and columns of the data are linearly dependent.
We provide limit results for the weak \con\ of these \pp es to Poisson or cluster Poisson processes. Based on
this \con\ we can also derive the limit laws of various \fct als of the ordered eigenvalues such as the
joint \con\ of a finite number of the largest order statistics, the joint limit law of the largest eigenvalue and the trace,
limit laws for successive ratios of ordered eigenvalues,
etc. We also develop some
limit theory for the singular values of the sample autocovariance matrices and their sums of squares. The theory is illustrated
for simulated data and for the components of the S\&P 500 stock index.
\end{abstract}

\maketitle
%\tableofcontents

\section{Estimation of the largest eigenvalues: an overview in the iid case}\label{sec:motivation}\setcounter{equation}{0}
\subsection{The light-tailed case}
%In recent years, \asy\ theory for large random matrices has a
%attracted a
%great deal of attention; see for example the monographs
%Bai and Silverstein \cite{bai:silverstein:2010} and Anderson et al.
%\cite{anderson:guionnet:zeitouni:2008}. The literature about {\em
%  heavy-tailed random matrices} is rather sparse. Here and in what
%follows, we call a random matrix {\em heavy-tailed} if its entries have
%\ds s with \regvary\ tails, typically with tail index below 4.
%The eigenvalues of heavy-tailed matrices
%with iid entries and dimension $p\times n$, where $p=p_n\to\infty$ and $p/n\to
%\gamma\in (0,\infty)$, were studied by Soshnikov \cite{soshnikov:2004,
%  soshnikov:2006} and  Auffinger et
%al. \cite{auffinger:arous:peche:2009}. In the latter reference and in
%Belinschi et al. \cite{belinschi:dembo:guionnet:2009},
%the sample covariance matrices of iid heavy-tailed \seq s were also studied.
%Bose et al. \cite{bose:hazra:saha:2009,bose:hazra:saha:2010}
%investigated the spectral norm of circulant-type heavy-tailed matrices
%and Ben Arous and Guionnet \cite{arous:guionnet:2008}
%studied the limits of heavy-tailed random Wigner matrices with
%infinite variance.
One of the exciting new areas of statistics is concerned with analyses of large data sets.
For such data one often studies the dependence structure via covariances and correlations.
In this paper we focus on one aspect: the estimation of the eigenvalues of the covariance matrix of a multivariate  \ts\
when the dimension $p$ of the series increases with the sample size $n$. In particular, we are interested in limit theory for the largest eigenvalues of the sample covariance matrix. This theory is closely related to topics from classical \evt\ such as
maximum domains of attraction with the corresponding normalizing and centering constants for maxima; cf. Embrechts et al.~\cite{embrechts:kluppelberg:mikosch:1997}, Resnick~\cite{resnick:2007,resnick:1987}.
Moreover, \pp\ \con\ with limiting Poisson and cluster Poisson processes enters in a natural way when one describes
the joint \con\ of the largest eigenvalues of the sample covariance matrix. Large deviation techniques find applications,
linking \evt\ with random walk theory and \pp\ \con . The objective of this paper is to illustrate some of the main
developments in random matrix theory for the particular case of the sample covariance matrix of multivariate \ts\ with
independent or dependent entries. We give special emphasis to the heavy-tailed case when \evt\ enters in a rather
straightforward way.
%We do not cover other parts of random matrix theory such as limit theory for the eigen-spectrum of random matrices with
%iid entries or other special types of matrices. Also in this case one can find many parallels with \evt.
\par
Classical multivariate \tsa\ deals with
observations which assume values in a $p$-dimensional space
where $p$ is ``relatively small'' compared to the sample size $n$.  With
the availability of large data sets $p$ can be ``large''
relative to $n$. One of the possible con\seq s is that standard asymptotics (such as the \clt )
break down and may even cause misleading results.
\par
The dependence
structure in multivariate data is often summarized by the covariance matrix which is typically
estimated by its sample analog.  For example, principal component analysis (PCA)
extracts principal component vectors corresponding to the largest
eigenvalues of the sample covariance matrix. The magnitudes of these eigenvalues provide an empirical \ms\ of the importance
of these components.
\par
If $p,n$ are fixed, a column of the $p\times n$ data matrix
\beao
\X =\X_n = \big(X_{it}\big)_{i=1,\ldots,p;t=1,\ldots,n}\,
\eeao
represents an observation of a $p$-dimensional \ts\ model with unknown parameters.
In this section we assume that the real-valued entries $X_{it}$ are iid, unless mentioned otherwise, and we write $X$ for a generic element.
One challenge is to infer information about the parameters
from the eigenvalues $\lambda_1,\ldots,\lambda_p$
of the {\em sample covariance matrix} $\X\X'$. In the notation we suppress the dependence of $(\la_i)$ on $n$ and $p$.
If $p$ and $n$ are finite and the columns of $\X$ are iid and multivariate normal,
Muirhead \cite{muirhead} derived a (rather complicated) formula for the joint distribution of the eigenvalues $(\lambda_i)$.
%Therefore a theoretical solution is already available in the finite $n$ case.
%The expression, however, contains an integral over the orthogonal group $\mathbb{O}(p)$.
%\begin{figure}[htb!]
%  \centering
 % \includegraphics[scale=0.5]{eigen_normal_qqplot1.pdf}
 % \caption{Largest eigenvalues of $500$ simulated sample covariance matrices against standard normal quantiles.}
 % \label{fig:QQ}
%\end{figure}

For $p$ fixed and $\nto$, assuming
$\X$ has centered normal entries and a diagonal covariance matrix $\Sigma$,
Anderson~\cite{anderson:1963}  derived the joint \asy\ density of $(\lambda_1, \ldots, \lambda_p)$.
%If the largest eigenvalue of the underlying covariance matrix appears with multiplicity one, he
%showed that the largest eigenvalue of the sample covariave matrix
%is asymptotically normal.
We quote from Johnstone~\cite{johnstone:2001}:
``The classic paper by Anderson~\cite{anderson:1963} gives the limiting joint distribution of the roots, but the
marginal distribution of the largest eigenvalue is hard to extract even in the null case'' (i.e., when the covariance matrix $\Sigma$ is
proportional to the identity matrix).
%If the largest entry of $\Sigma$ exceeds its second largest, one obtains a central limit type theorem for the largest eigenvalue of the sample covariance matrix.
%In Figure~\ref{fig:QQ} we chose $n=10000, p=50$ and $\Sigma =\diag(1024,1,\ldots,1)$, and plot the properly normalized and
%standardized largest eigenvalue against standard normal quantiles.

It turns out that limit theory for the largest eigenvalues becomes ``easier''
when the dimension $p$ increases with $n$.
Over the last 15 years there has been increasing interest in the case when $p=p_n\to\infty$ as $\nto$. In most of
the literature (exceptions are El Karoui \cite{elkaroui:2003}, Davis et al. \cite{davis:mikosch:pfaffel:2015,davis:pfaffel:stelzer:2014} and Heiny and Mikosch
\cite{heiny:mikosch:2015:iid})
one assumes that $p$ and $n$ grow at the same rate:
\beam\label{eq:gamma}
\dfrac p n\to \gamma\qquad \mbox{for some $\gamma \in (0,\infty)$.}
\eeam
%One explanation for the commonly used condition \eqref{eq:gamma} is that much of the theory was developed for $n \times n$ {\em Wigner matrices}, which is a Hermitean matrix $H$ such that $\{ H_{ij}:i<j \}$ and $\{ H_{ii}:0\le i\le n \}$ are two sets of iid random variables. If $p =\gamma n$, it turns out that the limit theory for the eigenvalues of Wigner matrices and sample covariance matrices is closely connected by the convergence of the empirical spectral distributions to the semicircle and the Mar\v cenko--Pastur law, respectively; see Bai and Silverstein \cite{bai:silverstein:2010}, Chapters 1-3.

In random matrix theory,  the convergence of the {\em empirical spectral distributions} $(F_{n^{-1}\X \X'})$ of a sequence $(n^{-1}\X \X')$ of non-negative definite matrices is the principle object of study. The empirical spectral distribution $F_{n^{-1}\X \X'}$ is constructed from the eigenvalues via
\begin{equation*}
F_{n^{-1}\X \X'}(x)= \frac{1}{p}\; \# \{ 1\le j\le p : n^{-1} \lambda_j \le x \}, \quad x\in \R,\quad n\ge1.
\end{equation*}
In the literature convergence results for the sequence of empirical spectral distributions are established under the assumption that $p$ and $n$ grow at the same rate.
Suppose that the iid entries $Z_{it}$ have mean $0$ and variance $1$. If \eqref{eq:gamma} holds, then, with probability one, $(F_{n^{-1}\X \X'})$ converges to the celebrated Mar\v cenko--Pastur law with absolutely  continuous part given by the density,
\begin{eqnarray}\label{eq:MP}
f_\gamma(x) =
\left\{\begin{array}{cc}
\frac{1}{2\pi x\gamma} \sqrt{(b-x)(x-a)} \,, & \mbox{if } a\le x \le b, \\
 0 \,, & \mbox{otherwise,}
\end{array}\right.
\end{eqnarray}\noindent
where $a=(1-\sqrt{\gamma})^2$ and $b=(1+\sqrt{\gamma})^2$. The Mar\v cenko--Pastur law has a point mass $1-1/\gamma$ at the origin if $\gamma>1$, cf. Bai and Silverstein~\cite[Chapter~3]{bai:silverstein:2010}. The point mass at zero is intuitively explained by the fact that, with probability $1$, $\min(p,n)$ eigenvalues $\lambda_i$ are non-zero. When $n=(1/\gamma) \; p$ and $\gamma >1$ one sees that the proportion of non-zero eigenvalues of the sample covariance matrix is $1/\gamma$ while the proportion of zero eigenvalues is $1-1/\gamma$.

While the finite second moment is the central assumption to obtain the Mar\v cenko--Pastur law as the limiting spectral distribution, the finite fourth moment plays a crucial role when studying the largest eigenvalues
\beam\label{eq:order}
\la_{(1)}\ge \cdots \ge \la_{(p)}
\eeam
of $\X\X'$, where we suppress the dependence on $n$ in the notation.

\begin{figure}[htb!]
  \centering
  \subfigure[Standard normal entries]{
    \includegraphics[scale=0.4]{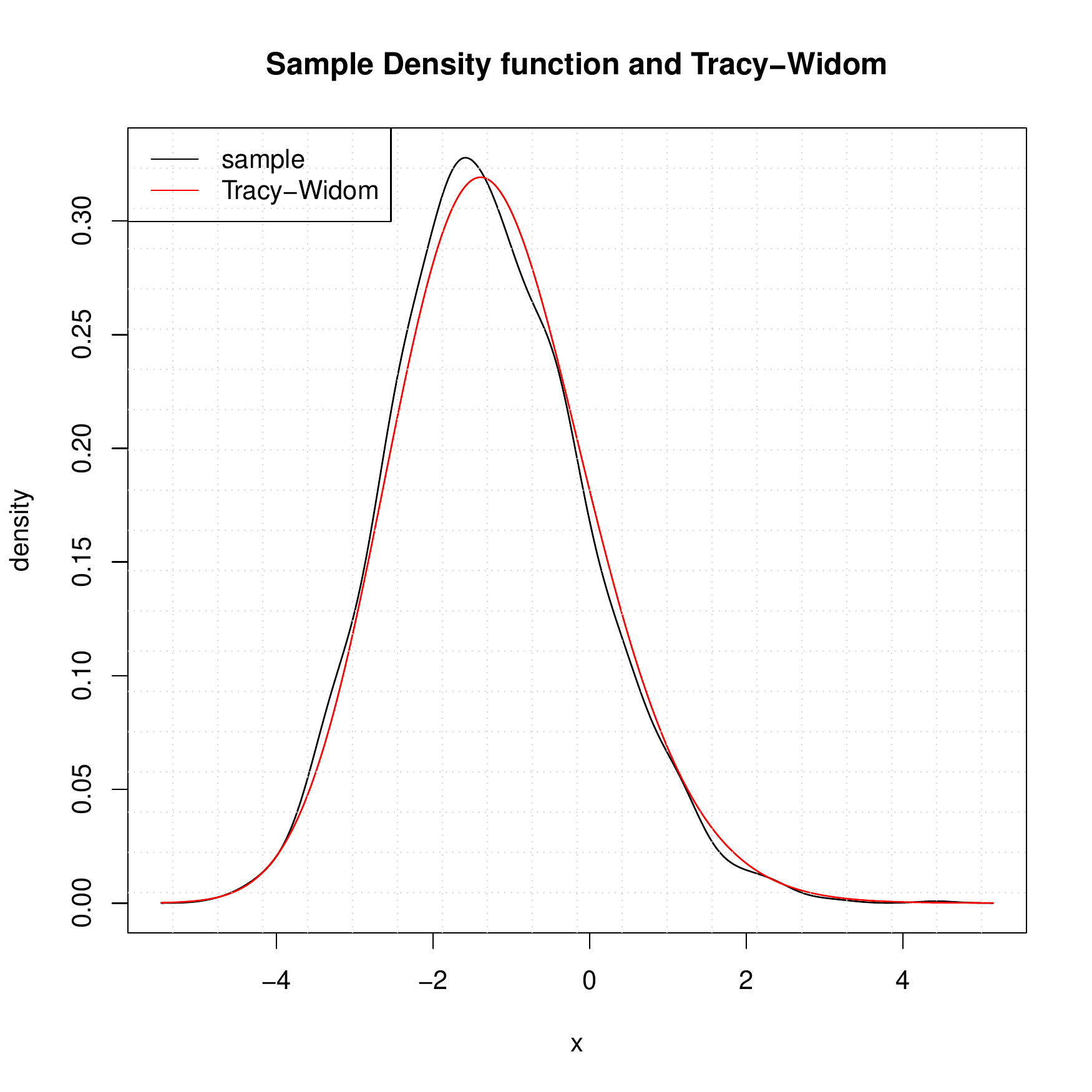}
  }
  \subfigure[Entry distribution: $\P(X=\sqrt{3}) = \P(X=-\sqrt{3} ) =
  1/6$, $\P(X=0)=2/3$. Note $\E X = 0$, { $\E[ X^2] = 1$, $\E [X^3] = 0$ and
  $\E[ X^4] = 3$, i.e., the first 4 moments of $X$ match those of the standard normal \ds .}]{
    \includegraphics[scale=0.4]{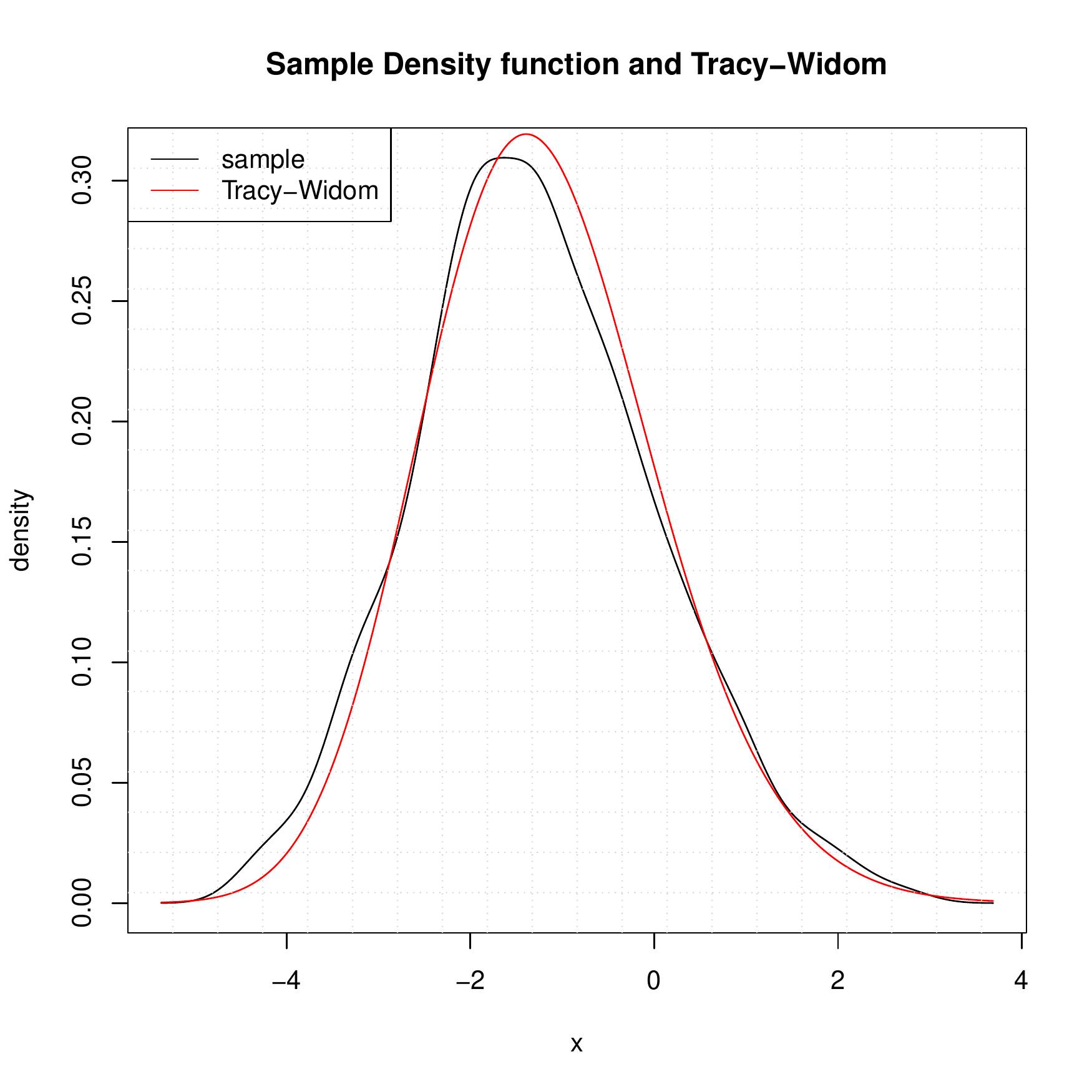}
  }
  \caption{Sample density function of the largest eigenvalue compared
    with the Tracy--Widom density function. The data matrix $ \X$ has
    dimension $200 \times 1000$. An ensemble of 2000 matrices is
    simulated.}
  \label{fig:normal-3point-TW}
\end{figure}

\begin{figure}[htb!]
  \centering
  \includegraphics[scale=0.5]{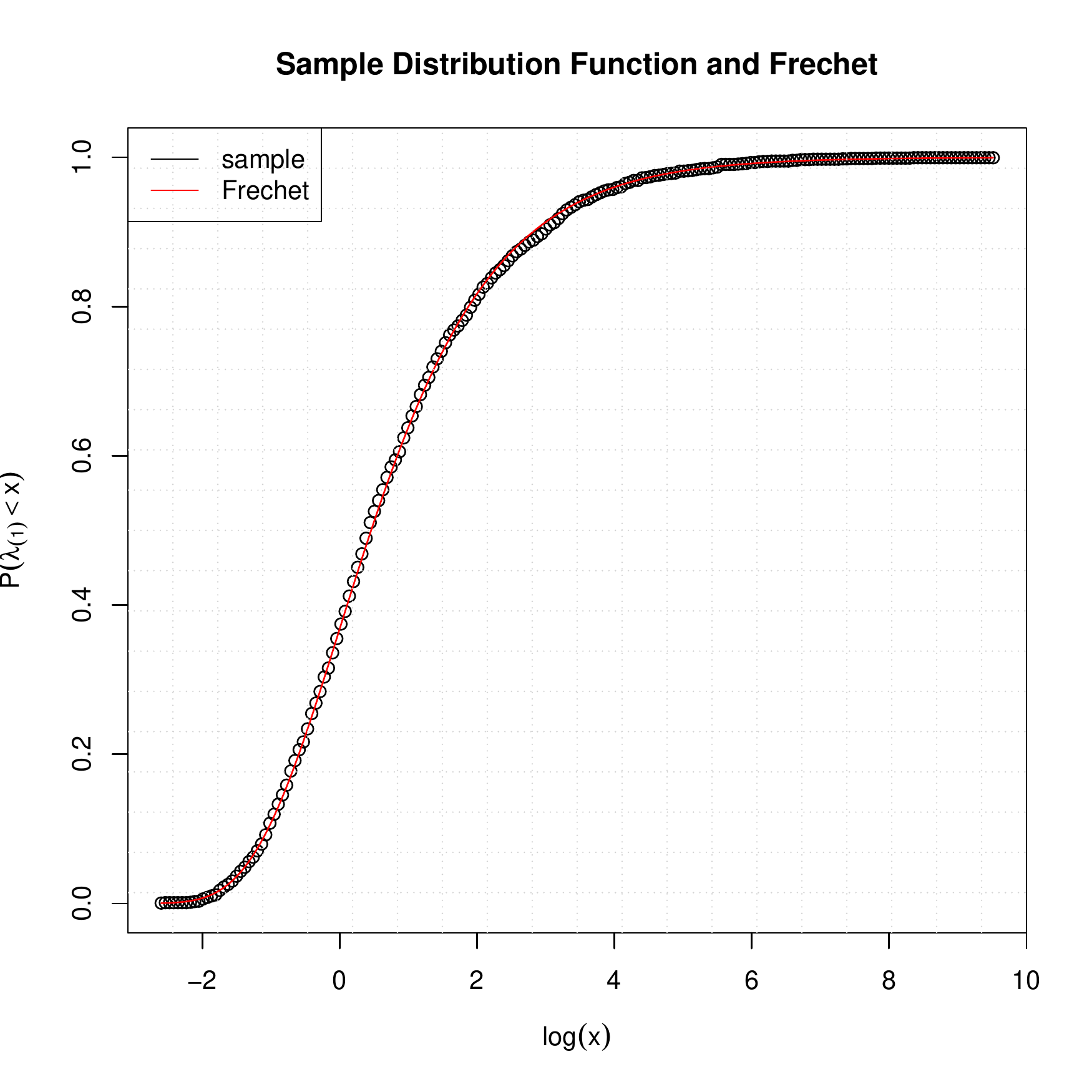}
  \caption{Sample distribution function of the largest eigenvalue $\la_{(1)}$
    compared to the \Frechet distribution (solid line) with $\alpha=1.6$. The data matrices have dimension
    $200 \times 1000$ and iid entries with infinite fourth moment. The results are based on 2000 replicates.}
  \label{fig:MyDist-Frechet}
\end{figure}

Assuming \eqref{eq:gamma} and iid entries $X_{it}$ with zero mean, unit variance and finite fourth moment,  Geman~\cite{geman}
showed that
\beam\label{eq:geman}
\dfrac {\la_{(1)}}{n} \stas \big(1+\sqrt{\gamma}\big)^2\,,\qquad \nto\,.
\eeam
Johnstone \cite{johnstone:2001} complemented this \slln\ by the corresponding \clt\ in the special case of iid standard normal
entries:
\beam\label{eq:tc}
 n^{2/3}\,\dfrac{(\sqrt{\gamma})^{1/3}}{\big(1+\sqrt{\gamma}\big)^{4/3}}\Big(\dfrac {\la_{(1)}}{n} -
\big(1+\sqrt{\tfrac pn }\big)^2\Big)
\std {\rm TW}\,,
\eeam
where the limiting \rv\ has a {\em Tracy--Widom \ds} of order 1. Notice that the centering
$\big(1+\sqrt{\tfrac pn }\big)^2$ can in general not be replaced by $(1+\sqrt{\gamma})^2$.
This \ds\ is ubiquitous in random matrix theory.
%It is defined via some ordinary differential equation; we refer to \cite{tracy:widom:2012} for a definition and properties.
Its distribution function $F_1$ is given by
\begin{equation*}
F_1(s) = \exp\Big\{
  -\frac{1}{2} \int_{s}^\infty [
    q(x) + (x - s) q^2(x)
 ] \dint x
\Big\}\,,
\end{equation*}
where $q(x)$ is the unique solution to the Painlev\'e II differential
equation
\begin{equation*}
  q''(x) = xq(x) + 2 q^3(x)\,,
\end{equation*}
where $ q(x)\sim {\rm Ai}(x)$ as $x \to \infty$ and Ai$(\cdot)$ is the Airy kernel; see Tracy and Widom~\cite{tracy:widom:2012} for details.
We notice that the rate $n^{2/3}$ compares favorably to the $\sqrt{n}$-rate in the classical \clt\ for sums
of iid finite variance \rv s.
The calculation of the spectrum is facilitated by the fact that the distribution of
the classical Gaussian matrix ensembles is invariant under orthogonal transformations. The corresponding
computation for non-invariant matrices with non-Gaussian entries is more complicated and was a major challenge for several years; a first step was made by Johansson \cite{johansson}.
Johnstone's result was extended to matrices $\X$ with iid non-Gaussian entries
by Tao and Vu \cite[Theorem~1.16]{tao09b}. Assuming that the first four moments of the entry \ds\ match  those of
the standard normal \ds , they showed \eqref{eq:tc} by
employing {\em Lindeberg's replacement method}, i.e., the iid non-Gaussian entries are replaced
step-by-step by iid Gaussian ones.
This approach is well-known from summation theory for \seq s of iid \rv s. Tao and Vu's result is a consequence of the so-called {\em Four Moment Theorem}, which describes the insensitivity of the eigenvalues with respect to changes in the distribution of the entries. To some extent (modulo the strong moment matching conditions) it shows the universality of Johnstone's limit
result \eqref{eq:tc}. Later we will deal with entries with infinite fourth moment. In this case, the weak limit
for the normalized largest eigenvalue $\la_{(1)}$ is distinct from the Tracy--Widom \ds : the classical Fr\'echet extreme value
\ds\ appears.
In Figure~\ref{fig:normal-3point-TW} we illustrate how the Tracy--Widom approximation works for Gaussian and non-Gaussian entries of $\X$
and in Figure~\ref{fig:MyDist-Frechet} we also illustrate that this approach fails when $\E[X^4]=\infty$.

Figure \ref{fig:normal-3point-TW} compares the sample density function
of the properly normalized largest eigenvalue estimated from 2000 simulated sample covariance matrices $\X\X'$ ($n=1000, p=200$) with
the Tracy--Widom density. If  $X$ has infinite fourth moment and further regularity conditions on the tail hold then the
Tracy--Widom limiting law needs to be replaced by the \Frechet
distribution; see Section~\ref{subsec:1.2} for details. Figure \ref{fig:MyDist-Frechet} illustrates this fact with a
simulated ensemble whose entries are distributed according to
the heavy-tailed distribution from \eqref{eq:distrsim} below with $\alpha = 1.6$.

\subsection{The heavy-tailed case}\label{subsec:1.2}
\begin{figure}[htb!]
\centering
\includegraphics[scale=0.5]{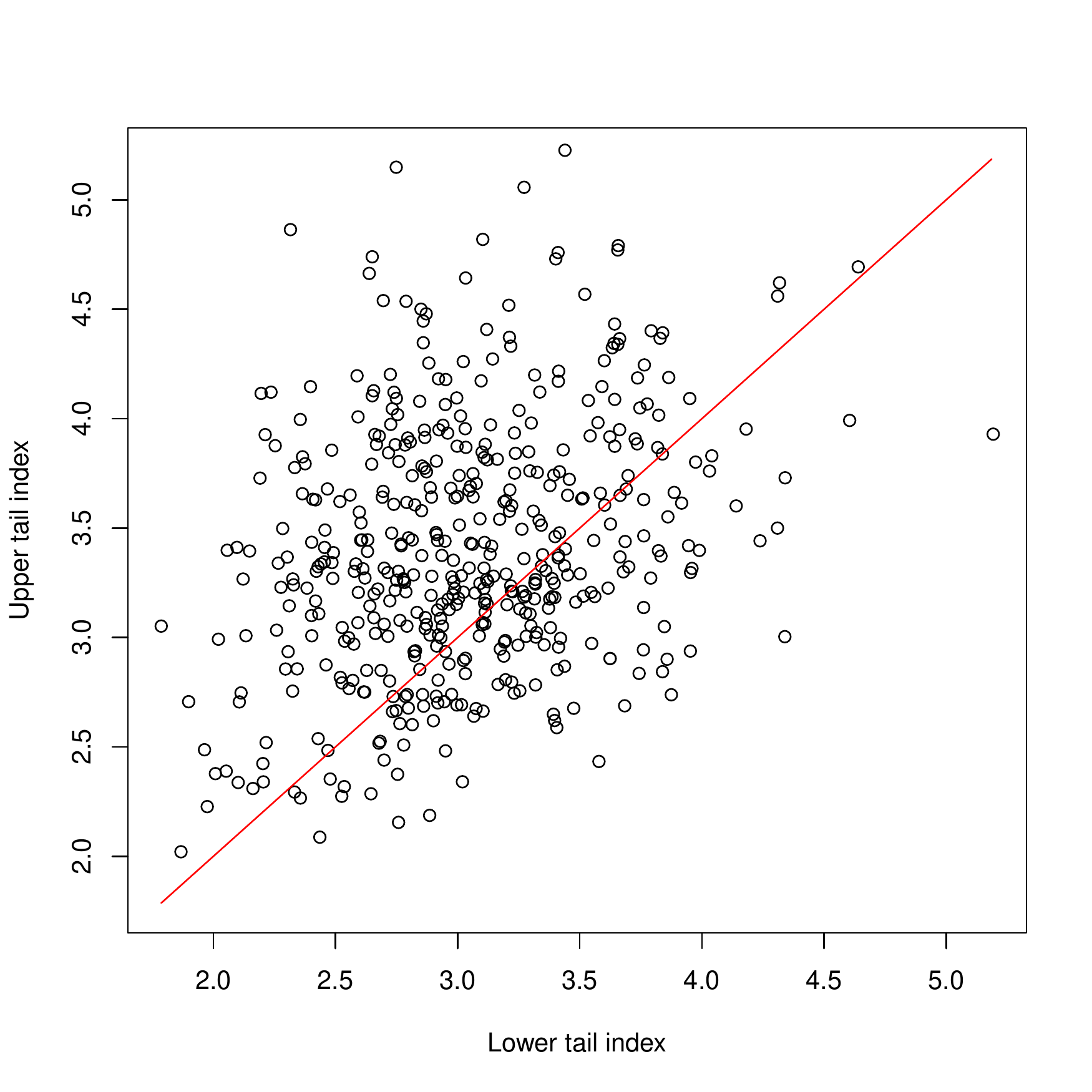}
\caption{Tail indices of log-returns of 478 \ts\ from the S\&P 500 index. The values  $(\hat \alpha_L,\hat \alpha_U)$ of the
lower and upper tail indices are provided by
Hill's estimator.  We also draw the line $\hat\alpha_U=\hat\alpha_L$.}
\label{fig:SP500_tail_indices}
\end{figure}

So far we focused on ``light-tailed'' $\X$
in the sense that its entries have finite fourth moment.  However, there is statistical evidence that the assumption
of finite fourth moment may be violated when dealing with data from
insurance, finance or telecommunications. We illustrate this fact
in Figure~\ref{fig:SP500_tail_indices} where we show the pairs $(\alpha_L,\alpha_U)$ of
lower and upper tail indices
of $p=478$  log-return series composing
the S\&P 500 index estimated from $n=1,345$ daily observations from 01/04/2010 to 02/28/2015.
This means we assume for every row  $(X_{it})_{t=1,\ldots,n}$ of $\X$ that the tails behave like
\beao
\P(X_{it}>x)\sim c_U\,x^{-\alpha_U}\qquad\mbox{and}\qquad \P(X_{it}<-x)\sim c_L\,x^{-\alpha_L}\,,\qquad \xto\,,
\eeao
for non-negative constants $c_L,c_U$. We apply the Hill estimator (see Embrechts et al.~\cite{embrechts:kluppelberg:mikosch:1997}, p.~330,
de Haan and Ferreira \cite{dehaan:ferreira:2006}, p.~69)
to the \ts\ of the gains and losses in a naive way,
neglecting the dependence and non-stationarity in the data; we also omit confidence bands.
From the figure it is evident that the majority of the return series have tail indices below
four, corresponding to an infinite fourth moment.
The behavior of the largest eigenvalue $\la_{(1)}$
changes dramatically when $\X$ has infinite fourth moment.
Bai and Silverstein~\cite{baisilv} proved for an $n\times n$ matrix $\X$ with iid centered entries
that
\beam\label{eq:wdfr}
\limsup_{\nto} \dfrac{\la_{(1)}}{n}=\infty \qquad {\rm a.s.}
\eeam
This is in stark contrast to Geman's result \eqref{eq:geman}.
%It is also unlikely that the
%Tracy--Widom \ds\ yields a good approximation to the \ds\ of the largest eigenvalue of the sample covariance matrix;
%see Figure~\ref{fig:MyDist-Frechet} for some simulation evidence.
\par
In the heavy-tailed case it is common to assume a {\em \regvar\ condition:}
\beam\label{eq:regvar}
\P(X>x)\sim p_+\,\dfrac{L(x)}{x^{\alpha}}\qquad\mbox{and}\qquad \P(X<-x) \sim p_-\,\dfrac{L(x)}{x^\alpha}\,,\qquad \xto\,,
\eeam
where $ p_\pm$ are non-negative constants \st\ $p_++p_-=1$ and $L$ is a \slvary\ \fct . In particular, if $\alpha<4$ we have $\E [X^4]= \infty$. The \regvar\ condition
on $X$ (we will also refer to $X$ as a \regvary\ \rv ) is needed for proving \asy\ theory for
the eigenvalues of $\X\X'$. This is similar to proving limit theory for sums
of iid \rv s with infinite variance stable limits; see for example Feller~\cite{feller}.
\par

In \eqref{eq:MP} we have seen that the sequence $(F_{n^{-1}\X \X'})$ of empirical spectral distributions converges to the Mar\v cenko--Pastur law if the centered iid entries possess a finite second moment. Now we will discuss the situation when the entries are still iid and centered, but have an infinite variance. Here we assume the entries to be regularly varying with index $\alpha \in (0,2)$.
Assuming \eqref{eq:gamma} with $\gamma\in (0,1]$ in this infinite variance case, Belinschi et al.~\cite[Theorem~1.10]{belinschi:dembo:guionnet:2009} showed that the sequence $(F_{a_{n+p}^{-2}\X \X'})$ converges with probability one to a non-random probability measure with density $\rho_{\alpha}^\gamma$ satisfying
\begin{equation*}
\rho_{\alpha}^\gamma(x) x^{1+\alpha/2} \to \frac{\alpha \gamma}{2(1+\gamma)}, \quad x \to \infty,
\end{equation*}
see also Ben Arous and Guionnet \cite[Theorem~1.6]{arous:guionnet:2008}.
The normalization $(a_k)$ is chosen \st\ $\P(|X|>a_k)\sim k^{-1}$ as $k\to\infty$. An application of the Potter bounds (see Bingham et al. \cite[p.~25]{bingham:goldie:teugels:1987}) shows that $a_{n+p}^2/n \to \infty$.

It is interesting to note that there is a phase change in the extreme eigenvalues in going from finite to infinite fourth moment, while the phase change occurs for the empirical spectral distribution going from finite to infinite variance.
%In the case of finite fourth moment, \eqref{eq:geman} and \eqref{eq:MP} imply that $\lambda_{(k)}/n$ has to converge to the right endpoint of the Mar\v cenko--Pastur law for any finite $k$. However, if we assume an infinite fourth moment but finite second moment the empirical spectral distributions  $(F_{{n}^{-1}\X \X'})$ still converge to the Mar\v cenko--Pastur law with the same support, while $\lambda_{(1)}/n$ will tend to infintity due to \eqref{eq:wdfr}. This means that the proportion of the normalized eigenvalues taking values larger than $(1+\sqrt{\gamma})^2$ has to be negligible; otherwise we would have a contradiction to the convergence of the empirical spectral distributions to the Mar\v cenko--Pastur law.

%Finally, if the second moment is infinity, i.e.~regularly varying entries with index $\alpha \in (0,2)$, then the much larger values of the eigenvalues are visible in the limiting spectral distribution even when using a much stronger normalization \cite[Theorem~1.10]{belinschi:dembo:guionnet:2009}.
\par

The theory for the largest eigenvalues of sample covariance matrices with heavy tails is less developed than in the light-tailed case.
Pioneering work for $\la_{(1)}$ in the case of iid \regvary\ entries $X_{it}$ with index $\alpha\in (0,2)$
is due to Soshnikov~\cite{soshnikov:2004,soshnikov:2006}. He showed the \pp\ \con\
\beam\label{eq:nn}
N_n=\sum_{i=1}^p \vep_{a_{np}^{-2}\la_i} \std N=\sum_{i=1}^\infty \vep_{\Gamma_i^{-2/\alpha}}\,,\qquad \nto\,,
\eeam
under the growth condition \eqref{eq:gamma} on $(p_n)$.
%The normalization $(a_k)$ is chosen \st\ $\P(|X|>a_k)\sim k^{-1}$ as $k\to\infty$.
Here
\beam\label{eq:Gamma}
\Gamma_i=E_1+\cdots + E_i\,,\qquad i\ge 1\,,
\eeam
and $(E_i)$ is an iid standard exponential \seq . In other words, $N$ is a Poisson point process on $(0,\infty)$ with mean \ms\
$\mu(x,\infty)= x^{-\alpha/2}$, $x>0$. Convergence in \ds\ of \pp es is understood in the sense of weak \con\
in the space of point \ms s equipped with the vague topology; see Resnick \cite{resnick:2007,resnick:1987}.
We can easily derive the limiting distribution of $a_{np}^{-2} \lambda_{(k)}$ for fixed $k\ge 1$ from \eqref{eq:nn}:
\begin{equation*}
\begin{split}
\lim_{\nto}\P(a_{np}^{-2} \lambda_{(k)}\le x)&= \lim_{\nto}\P(N_n(x,\infty)<k)
=  \P(N(x,\infty)<k) =\P(\Gamma_k^{-2/\alpha}\le x)\\ &= \sum_{s=0}^{k-1} \frac{\big(\mu(x,\infty)\big)^s}{s!} \e^{-\mu(x,\infty)}, \quad  x>0.
\end{split}
\end{equation*}
In particular,
\beao
\dfrac{\la_{(1)}}{a_{np}^2}\std \Gamma_1^{-\alpha/2}\,,\qquad \nto\,,
\eeao
where the limit has {\em Fr\'echet \ds } with parameter $\alpha/2$ and distribution function
\beao
\Phi_{\alpha/2}(x) =\ex^{-x^{-\alpha/2}}\,,\qquad x>0\,.
\eeao

\noindent We mention that the tail balance condition \eqref{eq:regvar} may be replaced in this case by the weaker assumption
$\P(|X|>x)= L(x) x^{-\alpha}$ for a \slvary\ \fct\ $L$. Indeed, it follows from the proof %(see also the comments below)
that only the squares $X_{it}^2$ contribute to the \pp\ limits of the eigenvalues $(\la_i)$. A con\seq\ of the
\cmt\ and \eqref{eq:nn} is the joint \con\ of the upper order statistics: for any $k\ge 1$,
\beao
a_{np}^{-2} \big(\la_{(1)},\ldots,\la_{(k)}\big)\std \big(\Gamma_1^{-2/\alpha},\ldots,\Gamma_k^{-2/\alpha}\big)\,,\qquad \nto\,.
\eeao

It follows from standard theory
for \pp es with iid points (e.g. Resnick \cite{resnick:2007,resnick:1987})
that \eqref{eq:nn} remains valid if
we replace  $N_n$  by the \pp\ $\sum_{i=1}^p\sum_{t=1}^n\vep_{X_{it}^2/a_{np}^2}$. Then we also have for any $k\ge 1$,
\beam\label{eq:nna}
a_{np}^{-2} \big(X_{(1),np}^2,\ldots,X_{(k),np}^2\big)\std \big(\Gamma_1^{-2/\alpha},\ldots,\Gamma_k^{-2/\alpha}\big)\,,\qquad \nto\,,
\eeam
where
\beao
X_{(1),np}^2\ge \cdots \ge X_{(np),np}^2
\eeao
denote the order statistics of $(X_{it}^2)_{i=1,\ldots,p;t=1,\ldots,n}$.
\par
Auffinger et al.~\cite{auffinger:arous:peche:2009} showed that \eqref{eq:nn}
remains valid under the \regvar\ condition \eqref{eq:regvar} for  $\alpha\in (2,4)$, the growth condition
\eqref{eq:gamma} on $(p_n)$ and the
additional assumption $\E [X]=0$. Of course, \eqref{eq:nna} remains valid as well.
Davis et al.~\cite{davis:pfaffel:stelzer:2014} extended these results to the case when the rows of $\X$
are iid linear processes with iid \regvary\ noise. The Poisson \pp\ \con\ result of \eqref{eq:nn} remains valid in this case.
Different limit processes can only be expected if there is dependence across rows and columns.
\par
In what follows, we refer to the  {\em heavy-tailed case} when we assume the \regvar\ condition \eqref{eq:regvar}
for some $\alpha\in (0,4)$.%{\blue This hypothesis is confirmed in}

\subsection{Overview}

The primary objective of this overview  is to make a connection between extreme value theory and the behavior of the largest eigenvalues of sample covariance matrices from heavy-tailed multivariate time series.  For time series that are linearly dependent through time and across rows, it turns out that the extreme eigenvalues are essentially determined by the extreme order statistics from an array of iid random variables.  The asymptotic behavior of the extreme eigenvalues is then derived routinely from classical extreme value theory.  As such, explicit joint distributions of the extreme order statistics can be given which yield a plethora of ancillary results.  
%These include the limit behavior of ratios of consecutive extreme eigenvalues, ratios of extreme order statistics to the sample variance, etc..
%%The underlying multivariate \ts\ will assumed to be heavy-tailed and linearly dependent through time and across the rows.We consider asymptotic theory for the largest eigenvalues of a sample covariance matrix,
%and we will also address the largest singular values of the sample autocovariance matrix.
%To our knowledge, the present paper is the first to consider bonafide dependence among the components in the time series which renders a multivariate analysis, such as PCA, meaningful.
Convergence of the point process of extreme eigenvalues, properly normalized, plays a central role in establishing the main results.
\par
In Section~\ref{sec:2} we continue the study of the case when the data matrix $\X$ consists of iid heavy-tailed entries.  We will consider power-law growth rates on the dimension $(p_n)$ that is more general than prescribed by \eqref{eq:gamma}.
In Section~\ref{sec:model} we introduce a model for $X_{it}$ which allows for linear dependence across the rows and through time.
The point process convergence of normalized eigenvalues is presented in Section~\ref{sec:mainresult}.
This result lays the foundation for new insight into the spectral behavior of the sample covariance matrix, which is the content of Section~\ref{sec:samplecov}.

Sections~\ref{sec:samplecov} and \ref{sec:possemidef} are devoted to {\em sample autocovariance matrices}.  Motivated by \cite{lam:yao}, we study the eigenvalues of sums of transformed matrices and illustrate the results in two examples. These results are applied to the time series of S\&P 500 in Section \ref{sec:sp500}.
%The technical proofs are collected in Section~\ref{sec:proof}. 
Appendix~\ref{appendix:A} contains useful facts about regular variation and point processes.

\section{General growth rates for $p_n$ in the iid heavy-tailed case}\label{sec:2}\setcounter{equation}{0}
This section is based on ideas in Heiny and Mikosch \cite{heiny:mikosch:2015:iid} where one can also find detailed proofs.
\subsection*{Growth conditions on $(p_n)$}\label{subsec:pn}
In many applications it is not realistic to assume
that the dimension $p$ of the data and the sample size $n$ grow at the same rate.
The aforementioned results of Soshnikov~\cite{soshnikov:2004,soshnikov:2006} and Auffinger et al. \cite{auffinger:arous:peche:2009} already indicate that the value $\gamma$ in the growth
rate~\eqref{eq:gamma} does not appear in the \ds al limits.
This obervation is in contrast to the light-tailed case; see \eqref{eq:geman}
and \eqref{eq:tc}.  Davis et al.~\cite{davis:mikosch:pfaffel:2015,davis:pfaffel:stelzer:2014}
and Heiny and Mikosch ~\cite{heiny:mikosch:2015:iid} allowed for more general rates for
$p_n\to\infty$ than linear growth in $n$.
Recall that $p=p_n\to\infty$ is the number of rows in the
matrix $\X_n$. We need to specify the growth rate of $(p_n)$ to ensure
a non-degenerate limit distribution of the normalized singular values of the sample autocovariance
matrices. To be precise, we assume
\begin{equation}\label{eq:p}
p=p_n=n^\beta \ell(n), \quad n\ge1,\tag{$C_p(\beta)$}
\end{equation}
where $\ell$ is a slowly varying function and $\beta\ge 0$. If $\beta =0$, we also assume $\ell(n) \to \infty$.
Condition~\ref{eq:p} is more general than the growth conditions in the literature; see
\cite{auffinger:arous:peche:2009,davis:mikosch:pfaffel:2015,davis:pfaffel:stelzer:2014}.
\begin{theorem}\label{thm:intro}
Assume that $\X=\X_n$ has iid entries satisfying the \regvar\ condition \eqref{eq:regvar} for some
$\alpha \in (0,4)$. If $\E[|X|]<\infty$ we also suppose that $\E [X]=0$. Let $(p_n)$ be an integer sequence
satisfying \ref{eq:p} with $\beta\ge 0$. In addition, we require
\begin{equation}\label{Cbeta}
\min(\beta,\beta^{-1})\in (\alpha/2-1,1] \qquad \mbox{ for } \alpha\in [2,4), \tag{$\widetilde{C}_\beta(\alpha)$}
\end{equation}
Then
\beam\label{eq:2}
\sum_{i=1}^p \vep_{a_{np}^{-2}\la_i}\std \sum_{i=1}^\infty \vep_{\Gamma_i^{-2/\alpha}}\,,\qquad \nto\,,
\eeam
where the \con\ holds in the space of point \ms s with state space $(0,\infty)$ equipped with the vague toplogy; see Resnick \cite{resnick:2007}.
\end{theorem}

\subsubsection*{A discussion of the case $\beta\in [0,1]$}
We mentioned earlier that in the heavy-tailed case, limit theory for the
largest eigenvalues of the sample covariance matrix is rather insensitive to
the growth rate of $(p_n)$ and that the limits are essentially
determined by the diagonal of this matrix. This is confirmed by the following result.

\begin{proposition}\label{prop:offdiagonal}
Assume that $\X=\X_n$ has iid entries satisfying the \regvar\ condition \eqref{eq:regvar} for some
$\alpha \in (0,4)$. If $\E[|X|]<\infty$ we also suppose that $\E [X]=0$. Then for any \seq\ $(p_n)$
satisfying \ref{eq:p} with $\beta\in [0,1]$ we have
\beao
a_{np}^{-2} \twonorm{\X \X' - \diag(\X \X')}\stp 0\,,\qquad\nto\,,
\eeao
where $\| \cdot\|_2$ denotes the spectral norm; see \eqref{specnorm} for its definition.
\end{proposition}

\par
Proposition~\ref{prop:offdiagonal} is not unexpected for two reasons:
\begin{itemize}
\item
It is well-known from classical theory (see Embrechts and Veraverbeke \cite{embrechts:veraverbeke:1982}) that for any iid \regvary\ non-negative \rv s
$Y,Y'$ with index  $\alpha'>0$, $Y\,Y'$ is \regvary\ with index $\alpha'$ while $Y^2$ is \regvary\ with index $\alpha'/2$.
Therefore $X^2$ and $X_{11}X_{12}$  are \regvary\ with indices $\alpha/2$ and $\alpha$, respectively.
\item
The aforementioned tail behavior is inherited by the entries of $\X\X'$ in the following sense.
By virtue of Nagaev-type \ld\ results for an iid \regvary\ \seq\ $(Y_i)$ with index $\alpha'\in (0,2)$ where we also assume that
$\E [Y_0]=0$ if $\E[|Y_0|]<\infty$
(see Theorem~\ref{thm:nagaev}) %S.V. Nagaev \cite{nagaev:1979} and
%Cline and Hsing \cite{cline:hsing:1998})
we have that $\P(Y_1+\cdots +Y_n >b_n)/(n \,\P(|Y_0|>b_n))$ converges to a non-negative constant
provided $b_n/a_n'\to\infty$, where $\P(|Y_0|>a_n')\sim n^{-1}$ as $\nto$. As a con\seq\ of the tail behaviors of $X_{it}^2$ and $X_{it}X_{jt}$ for $i\ne j$
and Nagaev's results we have for $(b_n)$ \st\ $b_n/a_n^2\to\infty$,
\beao
\dfrac{\P\big((\X\X')_{ij}> b_n\big)}{\P\big((\X\X')_{ii}- c_n> b_n\big)}\sim \dfrac{n\,\P( X_{11}X_{12}>b_n)}{n\,\P(X^2>b_n)}\to 0\,,\qquad \nto\,,
\eeao
where $c_n=0$ or $n\,\E[X^2]$ according as $\alpha\in (0,2)$ or $\alpha\in (2,4)$.
This means that the diagonal and off-diagonal entries of
$\X\X'$ inherit the tails of $X_{it}^2$ and $X_{it}X_{jt}$, $i\ne j$, respectively, above the high threshold $b_n$.
\end{itemize}
\par
Proposition~\ref{prop:offdiagonal} has some immediate con\seq s for the approximation of the eigenvalues
of $\X\X'$ by those of ${\rm diag}(\X\X')$. Indeed, let $C$ be a symmetric $p\times p$ matrix with
eigenvalues $\la_1(C),\ldots,\la_p(C)$ and
ordered eigenvalues
\beam \label{eq:weyl}
\la_{(1)}(C)\ge \cdots \ge \la_{(p)}(C)\,.
\eeam
Then for any symmetric $p\times p$ matrices $A,B$, by {\em Weyl's inequality} (see Bhatia \cite{bhatia:1997}),
\beao
\max_{i=1,\ldots,p}\big|\la_{(i)}(A+B)-\la_{(i)}(A)\big|\le \|B\|_2\,.
\eeao
If we now choose
$A+B=\X\X'$ and $A= \diag (\X\X')$ we obtain the following result.
\begin{corollary}\label{cor:687}
Under the conditions of Proposition~\ref{prop:offdiagonal},
\beao
a_{np}^{-2}\,\max_{i=1,\ldots,p}\big|\la_{(i)}-\la_{(i)}(\diag(\X\X'))\big|\stp 0\,,\quad\nto \,.
\eeao
\end{corollary}
Thus the problem of deriving limit theory for the order statistics of $\X\X'$ has been reduced to limit theory for the order
statistics of the iid row-sums
\beao
D_i^{\rightarrow}= (\X\X')_{ii}=\sum_{t=1}^n X_{it}^2\,,\qquad i=1,\ldots,p\,,
\eeao
which are the eigenvalues of $\diag(\X\X')$. This theory is completely described by the \pp es constructed from the points $D_i^\rightarrow/a_{np}^2$ $i=1,\ldots,p$. Necessary
and sufficient conditions for the weak \con\ of these \pp es are provided by  Lemma~\ref{lem:ppr}
which in combination with the Nagaev-type \ld\ results of  Theorem~\ref{thm:nagaev}
yield the following result; see also Davis et al.~\cite{davis:mikosch:pfaffel:2015}.
\begin{lemma}\label{lem:pp}
Assume the conditions of Proposition~\ref{prop:offdiagonal} hold. Then
\beao
\sum_{i=1}^p \vep_{a_{np}^{-2}(D_i^\rightarrow-c_n)}\std \sum_{i=1}^\infty \vep_{\Gamma_i^{-2/\alpha}}\,,\qquad \nto\,,
\eeao
where $(\Gamma_i)$ is defined in \eqref{eq:Gamma} and $c_n=0$ if $\E[D^\rightarrow]=\infty$ and $c_n=\E[D^\rightarrow]$ otherwise.
\end{lemma}
In this result, centering is only needed for $\alpha\in [2,4)$ when $n/a_{np}^2\not\to 0$. Under the additional condition \ref{Cbeta},
%\begin{equation}\label{Cbetaor}
%\beta\in (\alpha/2-1,1], \quad \mbox{ for } \alpha\in [2,4),
%\end{equation}
$n/a_{np}^2\to 0$ in view of the Potter bounds; see Bingham et al. \cite[p.~25]{bingham:goldie:teugels:1987}.
Combining Lemma~\ref{lem:pp} and Corollary~\ref{cor:687}, we conclude that Theorem~\ref{thm:intro} holds for $\beta\in [0,1]$.

%In  Davis et al. \cite{davis:mikosch:pfaffel;2015} relation \eqref{eq:2} was derived for $\alpha\in (0,1)$ even for \seq s
%$p_n\to\infty$ \st\ $p_n= O(\ex^{s_n/n})$, where $s_n=o(n)$ but with stronger limitations on the growth of $(p_n)$ for $[1,4]$.
%In the case $\alpha\in (2,4]$ they proved \eqref{eq:center} with $(\la_i)$ replaced by $(\la_i(\X\X'-n\E[ \X\X']))$. It is surprising that
%\eqref{eq:center} also holds since, in general, $\la_i(\X\X'-\E [\X\X'])\not = \la_i-n\,\E X$.

\subsubsection*{Extension to general $\beta$}%\label{subsec:generalb}
Next we explain that it suffices to consider only the case $\beta\in [0,1]$ and how to proceed when $\beta>1$.
The main reason is that
the $p \times p$ sample covariance matrix
$\X\X'$ and the $n \times n$ matrix $\X'\X$  have the same rank and their non-zero eigenvalues coincide; see Bhatia \cite[p.~64]{bhatia:1997}. When proving limit
theory for the eigenvalues  of the sample covariance matrix one may switch to $\X'\X$ and vice versa,
hereby interchanging the roles of $p$ and $n$. By switching to $\X'\X$, one basically replaces $\beta$ by $\beta^{-1}$. Since $\min(\beta,\beta^{-1})\in [0,1]$ for any $ \beta \ge 0$, one can assume without loss of generality that $\beta\in [0,1]$.  This trick allows one to extend results for
$(p_n)$ satisfying \ref{eq:p} with $\beta\in [0,1]$  to $\beta>1$. We illustrate this approach by providing the direct analogs of Proposition~\ref{prop:offdiagonal} and Corollary~\ref{cor:687}.

\begin{proposition}\label{prop:offdiagonal1}
Assume that $\X=\X_n$ has iid entries satisfying the \regvar\ condition \eqref{eq:regvar} for some
$\alpha \in (0,4)$. If $\E[|X|]<\infty$ we also suppose that $\E [X]=0$. Then for any \seq\ $(p_n)$
satisfying \ref{eq:p} with $\beta>1$ we have
\beao
a_{np}^{-2} \twonorm{\X' \X - \diag(\X' \X)}\stp 0\,,\qquad\nto\,,
\eeao
where $\| \cdot\|_2$ denotes the spectral norm.
\end{proposition}
Note that for $\beta>1$ we have $\lim_{\nto} p/n= \infty$. This means that
$\X'\X$ has  asymptotically a much smaller dimension than $\X\X'$ and therefore it is more convenient to work with $\X' \X$
when bounding the spectral norm.
\begin{corollary}\label{cor:6871}
Under the conditions of Proposition~\ref{prop:offdiagonal1},
\beao
a_{np}^{-2}\,\max_{i=1,\ldots,n}\big|\la_{(i)}-\la_{(i)}(\diag(\X'\X))\big|\stp 0\,,\quad\nto \,.
\eeao
\end{corollary}
Now, Theorem~\ref{thm:intro} for $\beta > 1$ is a consequence of Corollary~\ref{cor:6871}.
\par

%-----------------------------------------------------------------------
\section{Introducing dependence between the rows and columns}\label{sec:model}\setcounter{equation}{0}
For details on the results of this section, we refer to Davis et al.~\cite{davis:mikosch:pfaffel:2015}, Heiny and Mikosch \cite{heiny:mikosch:2015:iid} and Heiny et al.~\cite{heiny:mikosch:2016:noniid}. 

\subsection{The model}
When dealing with covariance matrices of a multivariate \ts\ $(\X_n)$ it is rather natural to assume dependence
between the entries $X_{it}$.
In this section we introduce a model which allows for {\em linear dependence}
between the rows and columns of $\X$:
\begin{equation}\label{eq:1}
X_{it}=\sum_{l\in \Z}\sum_{k\in \Z} h_{kl} Z_{i-k,t-l}\,,\qquad i,t\in\Z\,,
\end{equation}
where $(Z_{it})_{i,t\in \Z}$ is a field of iid \rv s and $(h_{kl})_{k,l\in\Z}$ is an array of real numbers.
Of course, linear dependence is restrictive in some sense. However, the particular dependence structure allows one to
determine those ingredients in the sample covariance matrix which contribute to its largest eigenvalues.
If the series in \eqref{eq:1} converges a.s. $(X_{it})$ constitutes a strictly stationary random field.
We denote generic elements of the $Z$- and $X$-fields by $Z$ and $X$, respectively. We assume that $Z$ is \regvary\ in the sense that
\begin{equation}\label{eq:27}
\P(Z>x)\sim p_+ \dfrac{L(x)}{x^{\alpha}}\quad\mbox{and}\quad  \P(Z\le -x)\sim p_-
\dfrac{L(x)}{x^{\alpha}}\,,\qquad \xto\,,
\end{equation}
for some tail index $\alpha>0$, constants $p_+,p_-\ge 0$ with $p_++p_-=1$ and a \slvary\ $L$. We will assume $\E[ Z]=0$ whenever $\E [Z^2]<\infty$.
Moreover, we require the summability condition
\begin{equation}\label{eq:2a}
\sum_{l \in \Z} \sum_{k\in \Z} |h_{kl}|^{\delta} <\infty
\end{equation}
for some $\delta\in (0,\min({\alpha/2},1))$ which ensures
the a.s.~absolute convergence of the series in \eqref{eq:1}. Under the conditions \eqref{eq:27} and \eqref{eq:2a}, the marginal and
\fidi s of the field $(X_{it})$ are \regvary\ with index $\alpha$; see
Embrechts et al. \cite{embrechts:kluppelberg:mikosch:1997}, Appendix A3.3. Therefore we also refer to $(X_{it})$ and $(Z_{it})$
as \regvary\ fields.
\par
The model \eqref{eq:1} was introduced by Davis et al. \cite{davis:pfaffel:stelzer:2014}, assuming the rows iid, and in the
present form by Davis et al. \cite{davis:mikosch:pfaffel:2015}.
\subsection{Sample covariance and autocovariance matrices}
From the field $(X_{it})$ we construct the $p\times n$ matrices
\begin{equation}\label{eq:26}
\X_n(s)= (X_{i,t+s})_{i=1,\ldots,p;t=1,\ldots,n}\,,\qquad s=0,1,2,\ldots\,,
\end{equation}
As before, we will write $\X=\X_n(0)$.
Now we can introduce the (non-normalized)
{\em sample autocovariance matrices}
\beam\label{eq:sample}
\X_n(0)\X_n(s)'\,,\qquad s=0,1,2,\ldots\,.
\eeam
We will refer to $s$ as the {\em lag}. For $s=0$, we obtain the {\em sample covariance matrix.}
In what follows, we will be interested in the asymptotic behavior (of \fct s) of the eigen- and singular values of the
sample covariance and autocovariance matrices in the heavy-tailed case. Recall that the {\em singular values} of a matrix $A$ are the square roots of the
eigenvalues of the non-negative definite matrix $AA'$ and its {\em spectral norm} $\twonorm{A}$ is its largest singular value. We notice that $\X_n(0)\X_n(s)'$ is not symmetric
and therefore its eigenvalues can be complex. To avoid this situation, we use
%\beao
%\X_n(0)\X_n(s)' +\X_n(s)\X_n(0)'
%\eeao
%(this matrix is symmetric and therefore it has real eigenvalues) or use
the squares
\beam\label{eq:squares}
\X_n(0)\X_n(s)' \X_n(s)\X_n(0)'
\eeam
%(this matrix is  non-negative definite  and therefore it has non-negative eigenvalues).
whose eigenvalues are the squares of the singular values of $\X_n(0)\X_n(s)'$.
The idea of using the sample autocovariance matrices and \fct s of their squares \eqref{eq:squares}
originates from a paper by Lam and Yao~\cite{lam:yao} who used a model different from \eqref{eq:1}.
This idea is quite natural in the context of \tsa .
\par
In Theorem~\ref{thm:mains} below, we provide a general approximation result for the ordered singular values of the sample
autocovariance matrices in the heavy-tailed case. This result is rather technical. To formulate it we introduce further notation.
As before, $p=p_n$ is any integer \seq\ converging to infinity.
\subsection{More notation}
Important roles are played by the quantities $(Z_{it}^2)_{i=1,\ldots,p;t=1,\ldots,n}$ and their order statistics
 \begin{equation}\label{eq:zorder}
Z_{(1),np}^2 \ge Z_{(2),np}^2 \ge  \ldots \ge Z_{(np),np}^2, \qquad n,p\ge 1\,.
\end{equation}
As important are the row-sums
\begin{equation}\label{eq:ll}
D_i^\rightarrow=D_i^{(n),\rightarrow}=\sum_{t=1}^n Z_{it}^2\,, \qquad
i=1,\ldots,p;\quad  n=1,2,\ldots\,,
\end{equation} with generic element $D^\rightarrow$ and their ordered values
\beam\label{eq:help6}
D_{(1)}^\rightarrow=D_{L_1}^\rightarrow\ge \cdots \ge D_{(p)}^\rightarrow=D_{L_p}^\rightarrow\,,
\eeam
where we assume without loss of generality that $(L_1,\ldots,L_p)$ is
a permutation of $(1,\ldots,p)$ for fixed~$n$.
\par
Finally, we introduce the column-sums
\beao
D_t^\downarrow=D_t^{(n),\downarrow}= \sum_{i=1}^p Z_{it}^2\,,\qquad t=1,\ldots,n;\; \quad p=1,2,\ldots\,,
\eeao
with generic element $D^\downarrow$ and we also adapt the notation from \eqref{eq:help6} to these quantities.
\subsubsection*{Matrix norms} For any $p\times n $ matrix $\A=(a_{ij})$, we will use the following norms:
\begin {itemize}
\item
{\em Spectral norm:}
\beam\label{specnorm}
\|\A\|_2=\sqrt{\la_{(1)}(\A\A')}\,,
\eeam
\item
{\em Frobenius norm:}
\beao
\frobnorm{\A}= \Big( \sum_{i=1}^p \sum_{j=1}^n |a_{ij}|^2\Big)^{1/2}\,.
\eeao
%\item
%{\em Max-row sum norm:}
%\beao
%\inftynorm{A}= \max_{i=1,\ldots,p} \sum_{j=1}^n |a_{ij}|\,.
%\eeao
%\item
%{\em Max-column sum norm:}
%\beao
%\|{A}\|_1 = \max_{j=1,\ldots,n} \sum_{i=1}^p |a_{ij}|\,.
%\eeao
\end{itemize}
We will frequently make use of the bound $\|\A\|_2\le \|\A\|_ F$. Standard references for matrix norms are \cite{belitski,bhatia:1997, horn, shores}.
\subsubsection*{Singular values of the sample autocovariance matrices} Fix integers $n\ge 1$ and $s\ge 0$. We recycle the $\la$-notation
for the singular values
$\la_1(s), \ldots, \la_p(s)$ of the sample autocovariance matrix $\X_n(0)\X_n(s)'$, suppressing the dependence on $n$.
Correspondingly, the  order statistics are denoted by
\begin{equation}\label{eq:sigma}
\la_{(1)}(s) \ge \cdots \ge \la_{(p)}(s)\,.
\end{equation}
When $s=0$ we typically write $\la_i$ instead of $\la_i(0)$.
\subsubsection*{The matrix $\M(s)$} We introduce some auxiliary matrices
derived from the coefficient matrix $\H=(h_{kl})_{k,l\in \Z}$:
\beao
\H(s)=(h_{k,l+s})_{k,l\in \Z}, \qquad \M(s)= \H(0)\H(s)'\qquad s\ge 0\,.
\eeao
Notice that
\begin{equation}\label{eq:m}
(\M(s))_{ij}= \sum_{l\in \Z} h_{i,l} h_{j,l+s}, \qquad i,j \in \Z .
\end{equation}
We denote the ordered singular values of $\M(s)$ by
\begin{equation}\label{eq:v1}
v_1(s) \ge v_2(s) \ge\cdots \,.
\end{equation}
Let $r(s)$ be the rank of $\M(s)$ so that $v_{r(s)}(s)>0$ while $v_{r(s)+1}(s)=0$ if $r(s)$ is
finite, otherwise $v_i(s)>0$ for all $i$. We also write $r=r(0)$.
\par
Under the summability condition \eqref{eq:2a} on $(h_{kl})$ for fixed $s\ge 0$,
\beam\label{eq:tracea}
\sum_{i=1}^\infty (v_i(s))^2 &=&  \frobnorm{\M(s)}^2= \sum_{i,j\in \Z} \sum_{l_1,l_2 \in \Z} h_{i,l_1} h_{j,l_1+s}h_{i,l_2}
h_{j,l_2+s}\nonumber\\
&\le& c \,\Big(\sum_{l_1,l_2 \in \Z} \sum_{i\in \Z} |h_{i,l_1} h_{i,l_2}|\Big)^2
\le c \,\sum_{l_1\in \Z} \sum_{i \in \Z} |h_{i,l_1}| <\infty\,.
\eeam
Therefore all singular values $v_i(s)$ are finite and the ordering
  \eqref{eq:v1}
is justified.
\par
{\em Here and in what follows, we write $c$ for any constant whose value is not of interest.}

\subsubsection*{Normalizing sequence}
We define $(a_k)$ by
\beao
\P(|Z|>a_k)\sim k^{-1}\,,\qquad k\to\infty\,,
\eeao
and choose the normalizing \seq\ for the singular values as $(a_{np}^2)$ for suitable \seq s $p=p_n\to\infty$.
\subsubsection*{Approximations to singular values}\label{subsec:defdelta}
We will give approximations to the singular values $\la_i(s)$  in
terms of the  $p$ largest ordered values for $s\ge 0$,
\beao
&&\delta_{(1)}(s)\ge \cdots \ge \delta_{(p)}(s)\,,\\
&&\gamma_{(1)}^\rightarrow(s)\ge \cdots \ge \gamma_{(p)}^\rightarrow(s)\,,\\
&&\gamma_{(1)}^\downarrow(s)\ge \cdots \ge \gamma_{(n)}^\downarrow(s)\,,
\eeao
from the sets
\beao
&&\big\{Z_{(i),np}^2 v_j(s)\,, i=1,\ldots,p\,;j=1,2,\ldots\big\}\,,\notag\\
&&\big\{D_{i}^\rightarrow v_j(s), i=1,\ldots,p\,;j=1,2,\ldots\big\}\,,\notag\\
&&\big\{D_{t}^\downarrow v_j(s), t=1,\ldots,n\,;j=1,2,\ldots\big\}\,,
\eeao
respectively.
%---------------------------------------------------------------------------
\subsection{Approximation of the singular values}\label{sec:mainresult}
%Furthermore, let $\tau_{(1)}(s)  \ge \ldots \ge \tau_{(p)}(s)$ be  the singular values of $\tfrac{1}{2}(\X_n(0)\X_n(s)'+\X_n(s)\X_n(0)' )$ and let $\tilde{\tau}_{(1)}(s) \ge  \ldots \ge\tilde{\tau}_{(p)}(s)$ be the singular values of
%\begin{equation*}
%\tfrac{1}{2}(\X_n(0)\X_n(s)'-\E[\X_n(0)\X_n(s)']+\X_n(s)\X_n(0)' -\E[\X_n(s)\X_n(0)'])
%\end{equation*}
In the following result we povide some useful approximations to the singular values of the sample autocovariance matrices of
the linear model \eqref{eq:1}.
\begin{theorem}\label{thm:mains}
Consider the linear process \eqref{eq:1} under
\begin{itemize}
\item
the \regvar\ condition \eqref{eq:27}
for some $\alpha\in (0,4)$,
\item the centering condition
$\E[Z]=0$ if $\E[|Z|]<\infty$,
\item
the summability condition
\eqref{eq:2a} on the coefficient matrix  $(h_{kl})$,
\item
the growth condition \ref{eq:p} on $(p_n)$ for some $\beta\ge 0$.
\end{itemize}
Then the following statements hold for $s\ge 0$:
\begin{enumerate}
\item We consider two disjoint cases:
$\alpha \in (0,2)$ and $\beta\in (0,\infty)$, or
$\alpha\in [2,4)$ and $\beta$ satisfying \ref{Cbeta}. Then
\begin{equation}\label{eq:mains1}
a_{np}^{-2} \max_{i=1,\ldots,p} |\la_{(i)}(s)-\delta_{(i)}(s)| \cip 0, \quad \nto.
\end{equation}
\item
Assume $\beta\in [0,1]$.
If $\alpha \in (0,2]$, $\E[Z^2]=\infty$ or $\alpha\in [2,4)$, $\E [Z^2]<\infty$ and $\beta\in (\alpha/2-1,1]$ then
\beao
a_{np}^{-2} \max_{i=1,\ldots,p} |\la_{(i)}(s)-\gamma_{(i)}^\rightarrow(s)| \cip 0, \quad \nto.
\eeao
Assume $\beta>1$. If $\alpha \in (0,2]$, $\E[Z^2]=\infty$ or $\alpha\in [2,4)$, $\E [Z^2]<\infty$ and $\beta^{-1}\in(\alpha/2-1,1]$. Then
\beao
a_{np}^{-2} \max_{i=1,\ldots,p} |\la_{(i)}(s)-\gamma_{(i)}^\downarrow(s)| \cip 0, \quad \nto.
\eeao
\end{enumerate}
\end{theorem}
\begin{remark}\em
The proof of Theorem~\ref{thm:mains} is given in Heiny et al.~\cite{heiny:mikosch:2016:noniid}. Part (2) of this  result
with more restrictive conditions on the growth rate of $(p_n)$ is contained in Davis et al.~\cite{davis:mikosch:pfaffel:2015}.
These proofs are very technical and lengthy.
\end{remark}
\begin{remark}\em
If we consider a random array $(h_{kl})$ independent of  $(X_{it})$ and
assume that the summability
condition \eqref{eq:2a} holds a.s.
then Theorem~\ref{thm:mains} remains valid conditionally on
$(h_{kl})$, hence  unconditionally in $\P$-probability; see also \cite{davis:mikosch:pfaffel:2015}.\end{remark}
\subsection{Point process \con }
Theorem~\ref{thm:mains}  and arguments similar to the proofs in Davis et al. \cite{davis:mikosch:pfaffel:2015}
enable one to derive the weak \con\ of
the point processes of the normalized singular values. Recall
the representation of the points $(\Gamma_i)$ of a unit rate homogeneous Poisson process on $(0,\infty)$
given in \eqref{eq:Gamma}. For $s\ge 0$, we define the point processes of the normalized singular values:
\begin{equation}\label{eq:ppdef}
N_n^{\lambda,s}=\sum_{i=1}^p \vep_{a_{np}^{-2}(\lambda_{(i)}(0),\ldots,\lambda_{(i)}(s))} \,.
\end{equation}

\begin{theorem}\label{cor:1}
Assume the conditions of Theorem~\ref{thm:mains}.
Then $(N_n^{\lambda,s})$ converge weakly in the space of point measures
with state space $(0,\infty)^{s+1}$ equipped with the vague topology.
If either $\alpha \in (0,2]$, $\E[Z^2]=\infty$ and $\beta \ge 0$,
or $\alpha \in [2,4)$, $\E[Z^2]<\infty$ and \ref{Cbeta} hold then
\begin{equation}\label{eq:pp}
N_n^{\lambda,s} \cid N= \sum_{i=1}^\infty
\sum_{j=1}^{\infty} \vep_{\Gamma_i^{-2/\alpha} (v_j(0),\ldots,v_j(s))}, \qquad \nto.
\end{equation}
\end{theorem}

\begin{proof}
Regular variation of $Z^2$ is equivalent to
\begin{equation}\label{eq:wwa}
n\,p\, \P(a_{np}^{-2} Z^2  \in \cdot ) \civ \mu(\cdot),
\end{equation}
where $\civ$ denotes vague convergence of Radon measures on  $(0,\infty)$ and the \ms\ $\mu$ is given by $\mu(x,\infty)=x^{-\alpha/2}$, $x>0$.
In view of Resnick \cite{resnick:1987}, Proposition~3.21, \eqref{eq:wwa} is equivalent
to the weak \con\ of the following \pp es:
\begin{equation*}
\sum_{i=1}^p\sum_{t=1}^n \vep_{a_{np}^{-2} Z_{it}^2}=\sum_{i=1}^{np} \vep_{a_{np}^{-2} Z^2_{(i),np}}\cid
\sum_{i=1}^\infty\vep_{\Gamma_i^{-2/\alpha} }=\widetilde{N}\,,\quad
\nto\,,
\end{equation*}
where the limit $\widetilde{N}$ is a Poisson random measure (PRM) with state space $(0,\infty)$ and
mean measure~$\mu$.
\par
Since $a_{np}^{-2} Z^2_{(p),np}\cip 0$ as $\nto$,
the point processes $\sum_{i=1}^{p} \vep_{a_{np}^{-2} Z^2_{(i),np}}$ converge weakly to the same PRM:
\begin{equation}\label{eq:ppc}
\sum_{i=1}^{p} \vep_{a_{np}^{-2} Z^2_{(i),np}}\cid
\sum_{i=1}^\infty\vep_{\Gamma_i^{-2/\alpha} }\,,\quad
\nto\,.
\end{equation}
A continuous mapping argument together with the fact that $\sum_{i=1}^\infty (v_i(s))^2<\infty$ (see  \eqref{eq:tracea})
shows that
\begin{equation*}
\sum_{j=1}^{ \infty}\sum_{i=1}^p\vep_{a_{np}^{-2} Z^2_{(i),np}(v_j(0),\ldots,v_j(s))}\cid
\sum_{j=1}^{ \infty}\sum_{i=1}^\infty\vep_{\Gamma_i^{-2/\alpha}  (v_j(0),\ldots,v_j(s))}\,.
\end{equation*}
If the assumptions of part (1) of Theorem~\ref{thm:mains} are satisfied
an application of \eqref{eq:mains1} (also recalling the definition of $(\delta_{(i)}(s))$) shows that \eqref{eq:ppc} remains
valid with the points  $(a_{np}^{-2} Z^2_{(i),np}(v_j(0),\ldots,v_j(s)))$ replaced by
$(a_{np}^{-2} (\lambda_{(i)}(0),\ldots,\lambda_{(i)}(s))$.
\par
The only cases which are not covered by Theorem~\ref{thm:mains}(1) are $\alpha\in (0,2)$, $\beta=0$ and $\alpha=2$, $\E[Z^2]=\infty$, $\beta\ge0$.
In these cases we get from Theorem~\ref{thm:nagaev} that
\begin{equation*}
 p \, \P(a_{np}^{-2} D^\rightarrow>x ) \sim p\,n\, \P( Z^2>a_{np}^{2} x )\to \mu(x,\infty)\,,\quad x>0\,,
\end{equation*}
i.e., $
 p \, \P(a_{np}^{-2} D^\rightarrow\in \cdot ) \stackrel{v}{\rightarrow} \mu(\cdot)$.
It follows from Lemma~\ref{lem:ppr} that
$\sum_{i=1}^p\vep_{a_{np}^{-2} D_i^\rightarrow}\std \widetilde{N}$. As before,
a continuous mapping argument in combination with
the approximation obtained in Theorem~\ref{thm:mains}(2)
justifies the replacement of the points $(a_{np}^{-2} D_{(i)}^\rightarrow(v_j(0),\ldots,v_j(s)))$ by
$(a_{np}^{-2} (\lambda_{(i)}(0),\ldots,\lambda_{(i)}(s)))$ in the case $\beta\in [0,1]$. If $\beta>1$ one has to work with the
quantities $(D_i^\downarrow)_{i=1,\ldots,n}$ instead of  $(D_i^\rightarrow)_{i=1,\ldots,p}$ and one may follow the same argument as above.
This finishes the proof.
\end{proof}

\section{Some applications}\setcounter{equation}{0}
\subsection{Sample covariance matrices}\label{sec:samplecov}%\setcounter{equation}{0}
The sample covariance matrix $\X_n(0)\X_n(0)'=\X\X'$ is a non-negative definite matrix. Therefore
its eigenvalues and singular values coincide.  Moreover, $v_j=v_j(0)$, $j\ge 1$, are the eigenvalues of $\M=\M(0)$.
\par
Theorem~\ref{thm:mains}(1) yields an approximation of the ordered eigenvalues $(\la_{(i)})$ of $\X\X'$ by
the quantities $(\delta_{(i)})$ which are derived from the order statistics of $(Z_{it}^2)$.
Part (2) of this result provides an approximation of $(\la_{(i)})$ by the quantities
$(\gamma_{(i)}^{\rightarrow/\downarrow})$ which are derived from the order statistics
of the partial sums $(D_i^{\rightarrow/\downarrow})$.
\par
In the following example we illustrate the
quality of the two approximations.
\begin{example}\label{ex:xiao}\em
We choose a Pareto-type \ds\ for $Z$ with density
\beam\label{eq:distrsim}
f_Z(x) =
\left\{\begin{array}{cc}
 \frac{\alpha}{(4|x|)^{\alpha + 1}}\,, & \mbox{if } |x| > 1/4 \\
1\,, & \mbox{otherwise.}
\end{array}\right.
\eeam
We simulated $20,000$ matrices $\X_n$ for $n=1,000$ and $p=200$ whose
iid entries have this density. We assume $\beta=1$.
Note that $\M=\M(0)$ has rank one and $v_1=1$.
The estimated densities of the deviations $a_{np}^{-2}(\lambda_{(1)}-D_{(1)}^\rightarrow)$ and
$a_{np}^{-2}(\lambda_{(1)}-Z^2_{(1),np})$  based on the simulations are shown in
Figure~\ref{fig:lambda_comparison}. The approximation error is very small indeed.
\begin{figure}[htb!]
  \centering
  \includegraphics[scale=0.40]{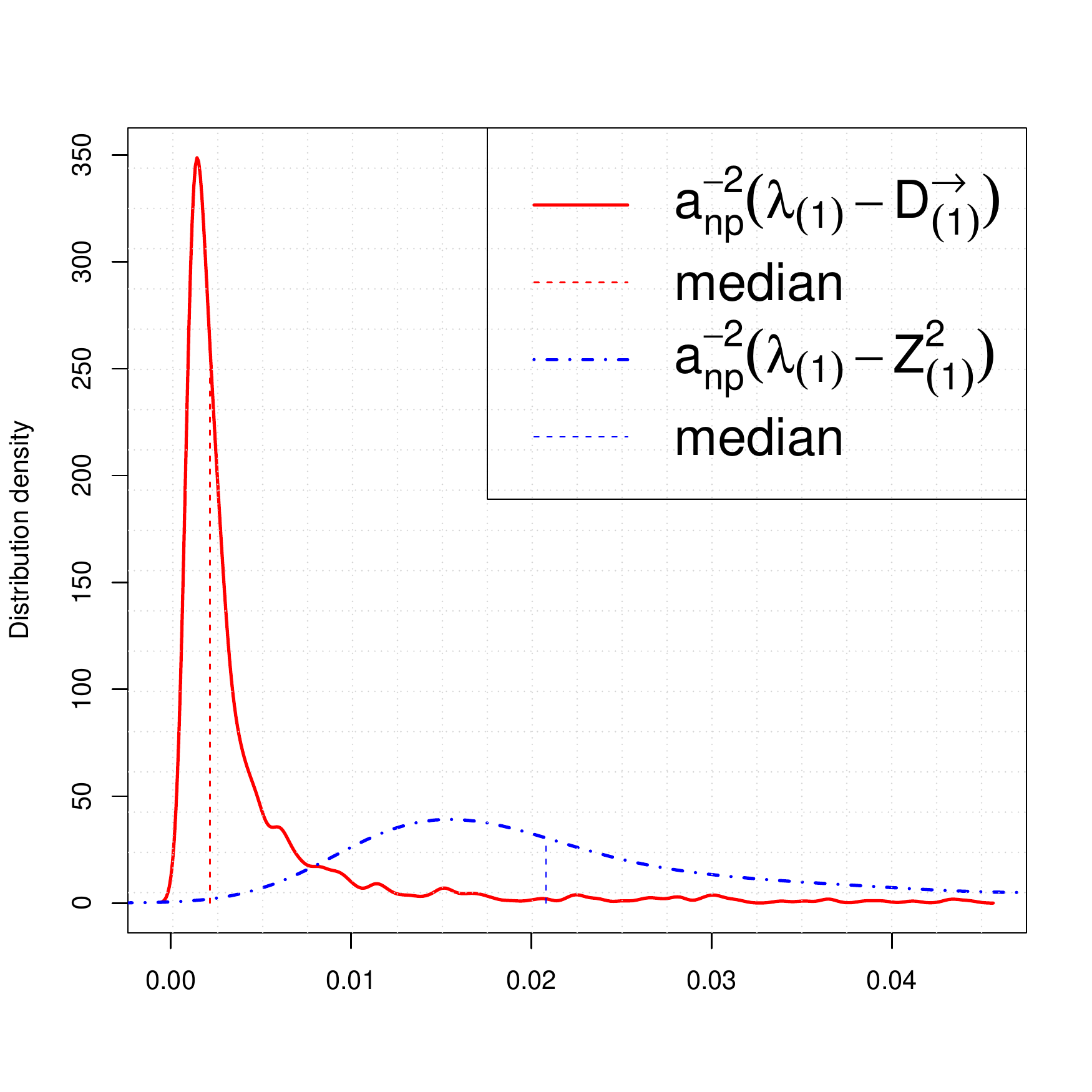}
  \caption{Density of the approximation errors for the eigenvalues of $a_{np}^{-2}\X\X'$.
The entries of $\X$ are iid with density \eqref{eq:distrsim} and $\alpha=1.6$.}
  \label{fig:lambda_comparison}
\end{figure}
According to the theory,
\beao
a_{np}^{-2}\sup_i|D_{(i)}^\rightarrow-\lambda_{(i)}| +a_{np}^{-2}\sup_i |Z^2_{(i),np}-\lambda_{(i)}|\stp 0\,,
\eeao
but  for finite $n$ the $(D_{(i)}^\rightarrow)$ sequence yields a better approximation to $(\la_{(i)})$. By construction,
the considered differences
have a tendency to be positive but Figure~\ref{fig:lambda_comparison} also shows that the median of the
approximation error for $a_{np}^{-2}(\lambda_{(1)}-D_{(1)}^\rightarrow)$ is almost zero.
\end{example}
\par
Theorem~\ref{cor:1} and the continuous mapping theorem immediately yield results about
the joint convergence of the largest eigenvalues of the matrices $a_{np}^{-2}\X_n\X_n'$ for $\alpha\in (0,2)$ and
$\alpha\in (2,4)$ when $\beta$ satisfies \ref{Cbeta}.
For fixed $k\ge 1$ one gets
\begin{equation*}
a_{np}^{-2} \big(\la_{(1)},\ldots,\la_{(k)}\big)\cid
\big(d_{(1)},\ldots, d_{(k)}\big)\,,
\end{equation*}
where $d_{(1)}\ge \cdots\ge d_{(k)}$ are the $k$ largest ordered
values of the set $\{\Gamma_i^{-2/\alpha} v_j, i=1,2,\ldots,j=1,\ldots,r\}$.
The continuous mapping theorem yields for $k\ge 1$,
\begin{equation}\label{eq:u}
\dfrac{\la_{(1)}}{\la_{(1)}+\cdots+ \la_{(k)}}\cid
\dfrac{d_{(1)}}{d_{(1)}+\cdots+ d_{(k)}}\,,\quad \nto\,.
\end{equation}

An application of the continuous mapping theorem to  the distributional convergence of the point processes in Theorem~\ref{cor:1} in the spirit of Resnick
\cite{resnick:2007}, Theorem 7.1, also yields the following result; see Davis et al. \cite{davis:mikosch:pfaffel:2015} for a proof and a similar result in the case $\alpha \in (2,4)$.
\begin{corollary}\label{cor:1q} Assume the conditions of Theorem~\ref{thm:mains}.
If $\alpha \in (0,2]$ and $\E[Z^2]=\infty$, then
\begin{equation*}
a_{np}^{-2}\Big(\la_{(1)},\sum_{i=1}^p \la_i\Big) \cid
\Big( v_1\,\Gamma_1^{-2/\alpha}\,,\sum_{j=1}^r v_j\,\sum_{i=1}^\infty \Gamma_i^{-2/\alpha}\Big)\,,
\end{equation*}
where $\Gamma_1^{-2/\alpha}$ is Fr\'echet $\Phi_{\alpha/2}$-distributed.
and $\sum_{i=1}^\infty \Gamma_i^{-2/\alpha}$ has the distribution
of a positive $\alpha/2$-stable random variable.
In particular,
\begin{equation}\label{eq:limit}
\dfrac{\la_{(1)}}{\la_{1}+\cdots+ \la_{p}}\cid
\dfrac{v_1}{\sum_{j=1}^r v_j}\;
\dfrac{\Gamma_1^{-2/\alpha}}{\sum_{i=1}^\infty\Gamma_i^{-2/\alpha}}\,,\quad
\nto\,.
\end{equation}
\end{corollary}

\begin{remark}\label{rem:4.5}\em
The ratio
\begin{equation*}
\dfrac{\la_{(1)}+\cdots +\la_{(k)}}{\la_1+\cdots+\la_p}, \qquad k\ge 1\,,
\end{equation*}
plays an important role in PCA. It reflects the proportion of the total variance in the data that we can explain by the first $k$ principal components.
It follows from Corollary~\ref{cor:1q} that for fixed $k\ge 1$,
\begin{equation*}
\dfrac{\la_{(1)}+\cdots \la_{(k)}}{\la_1+\cdots+\la_p}\cid
\dfrac{d_{(1)}+\cdots +d_{(k)}} {d_{(1)}+d_{(2)}+\cdots}\,.
\end{equation*}
Unfortunately, the limiting variable does in general not have a clean
form.
An exception is the case when $r=1$; see Example~\ref{exam:separable}.
Also notice that the trace of $\X\X'$ coincides with $\la_1+\cdots+\la_p$.
\end{remark}

To illustrate the theory we consider a simple moving average example
taken from Davis et al.~\cite{davis:mikosch:pfaffel:2015}.
\begin{example} \rm
\begin{figure}[htb!]
  \centering
  \subfigure[iid data]{
    \includegraphics[scale=0.4]{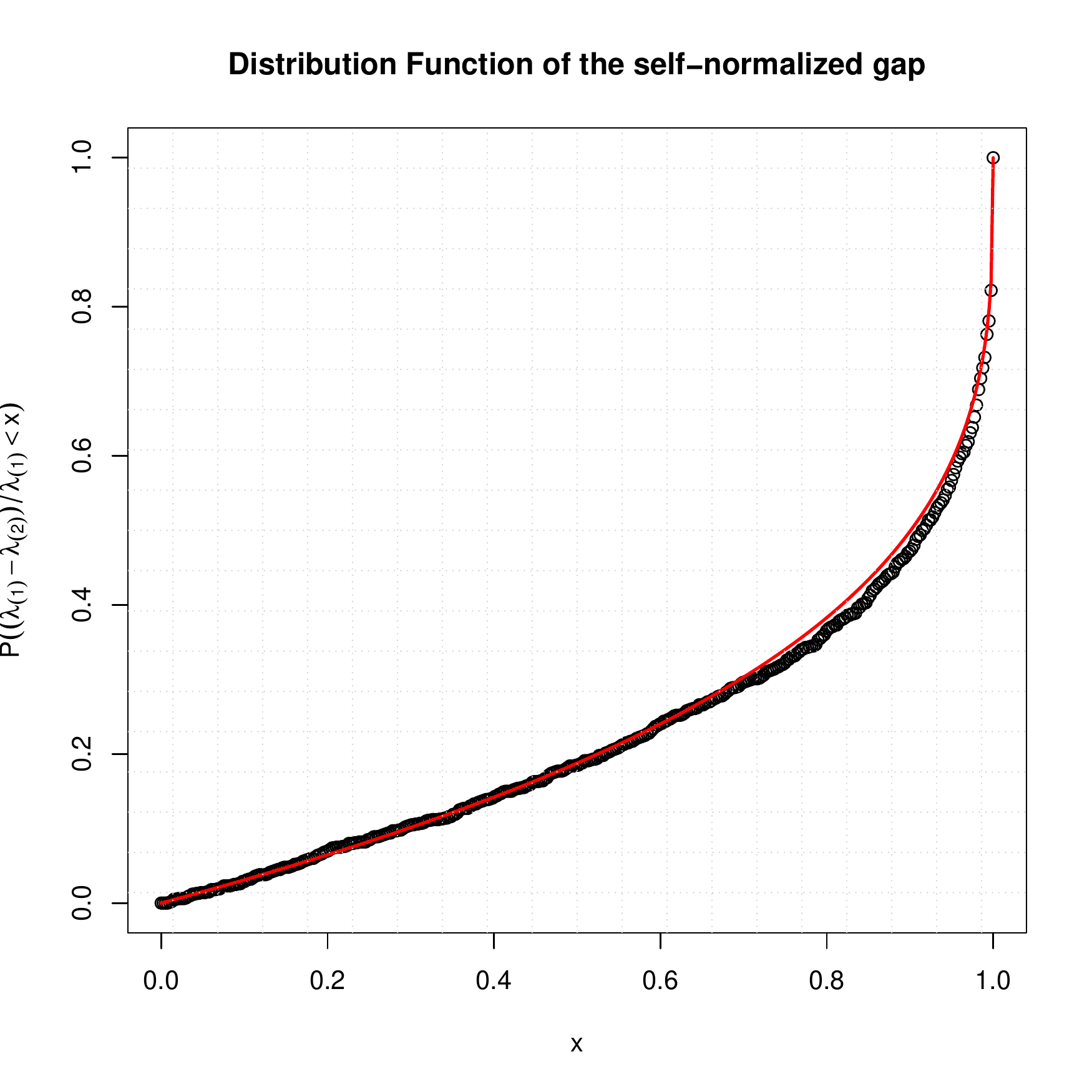}
  }
  \subfigure[data from model \eqref{ex:4.6}]{
    \includegraphics[scale=0.4]{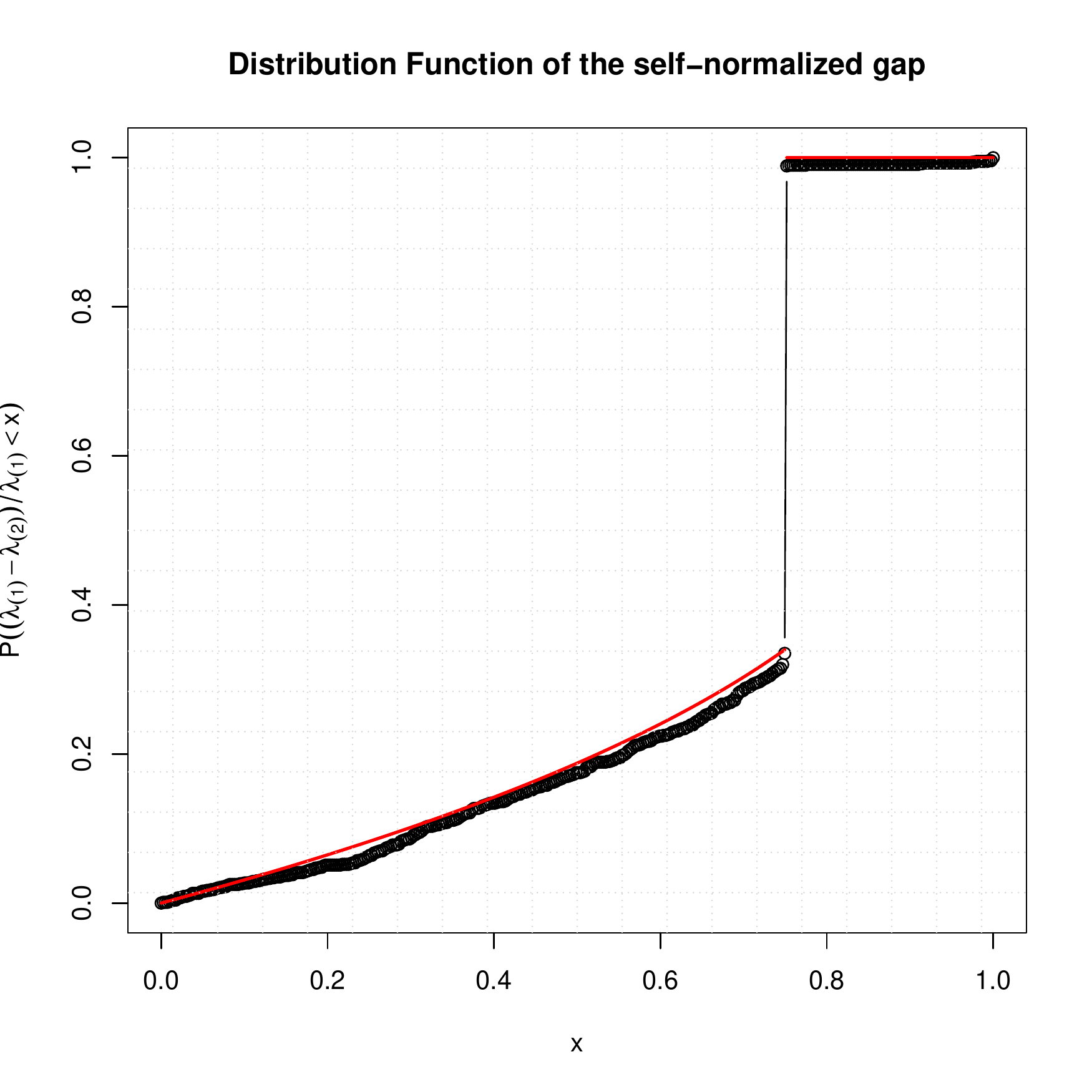}
  }
  \caption{Distribution function of $(\lambda_{(1)} - \lambda_{(2)})/\lambda_{(1)}$ for iid data (left) and
    data generated from the model \eqref{ex:4.6} (right). In each graph we compare the empirical \ds\ \fct\
(dotted line, based on
1000 simulations of $200 \times 1000$ matrices with $Z$-\ds\ \eqref{eq:distrsim}) with the theoretical curve (solid line).}
  \label{fig:ProbMass}
\end{figure}
Assume that $\alpha\in (0,2)$ and
\begin{equation}\label{ex:4.6}
X_{it}= Z_{it}+ Z_{i,t-1}-2 (Z_{i-1,t}- Z_{i-1,t-1})\,,\quad i,t\in\Z\,.
\end{equation}
In this case, the non-zero entries of $\H$ are
\begin{equation*}
h_{00}=1, h_{01}=1,h_{10}=-2 \quad \mbox{ and }\quad h_{11}=2.
\end{equation*}
Hence $\M=\H\H'$ has the positive eigenvalues
$v_1=8$ and $v_2=2$. The limit point process in \eqref{eq:pp} is
\begin{equation*}
N=\sum_{i=1}^\infty \vep_{8\Gamma_i^{-2/\alpha}}+ \sum_{i=1}^\infty \vep_{2\Gamma_i^{-2/\alpha}}\,,
\end{equation*}
so that
\begin{equation*}
a_{np}^{-2}\,\big(\la_{(1)},\la_{(2)}\big) \cid
\big(8\Gamma_{1}^{-2/\alpha}, 2\Gamma_1^{-2/\alpha}\vee 8
\Gamma_2^{-2/\alpha}\big)\,.
\end{equation*}
Using the fact that $U=\Gamma_1/\Gamma_2$ has a uniform \ds\ on $(0,1)$ we calculate
\begin{equation*}
\P(2\Gamma_1^{-2/\alpha}>8\Gamma_2^{-2/\alpha})= \P(\Gamma_1/\Gamma_2<2^{-\alpha})= 2^{-\alpha} \in (1/4,1).
\end{equation*}
In particular, we have for the normalized spectral gap
\begin{equation*}
a_{np}^{-2} \big(\la_{(1)}-\la_{(2)}\big)\cid
6 \,\Gamma_1^{-2/\alpha} \1_{\{\Gamma_1 4^{\alpha/2 }<\Gamma_2\}}+
8\,
\big(\Gamma_1^{-2/\alpha}-\Gamma_2^{-2/\alpha}\big)\1_{\{\Gamma_14^{\alpha/2
  }> \Gamma_2\}}
\end{equation*}
and for the self-normalized spectral gap (see also Example~\ref{exam:spectralgap} for a detailed analysis)
\begin{equation*}
\begin{split}
\dfrac{\la_{(1)}-\la_{(2)}}{\la_{(1)}}&\cid
\dfrac{6}{8} \, \1_{\{\Gamma_1 2^{\alpha }<\Gamma_2\}}+
\big(1-(\Gamma_1/\Gamma_2)^{2/\alpha}\big)\1_{\{\Gamma_12^{\alpha
  }> \Gamma_2\}}\\
& =\dfrac{3}{4}\,  \1_{\{U 2^{\alpha }<1\}}+
\big(1-U^{2/\alpha}\big)\1_{\{U2^{\alpha
  }> 1\}}=Y\,.
\end{split}
\end{equation*}
The limit  distribution of the spectral gap
has an atom at $3/4$ with probability $2^{-\alpha}$, i.e.,~$\P(Y=3/4)=2^{-\alpha}$, and
\begin{equation*}
\P(Y\le x)= 1- (1-x)^{\alpha/2},\qquad x\in (0,3/4).
\end{equation*}
In the iid case the limit distribution of the self-normalized spectral gap has distribution function $F(x) = 1- (1-x)^{\alpha/2}$ for $x\in [0,1]$. This means that the atom disappears if the entries are iid. Figure~\ref{fig:ProbMass} compares the distribution function of $Y$
with $F$ for $\alpha=0.6$; the atom at $3/4$ is clearly visible.

Along the same lines, we also have
\begin{equation*}
(a_{np}^{-2}\lambda_{(1)},\lambda_{(2)}/\lambda_{(1)}) \cid(
8\,\Gamma_1^{-2/\alpha},\frac{1}{4}\,\1_{\{U<2^{-\alpha}\}}+U^{2/\alpha}\,\1_{\{U\ge
  2^{-\alpha}\}})\,
\end{equation*}
and hence the limit distribution of $\lambda_{(2)}/\lambda_{(1)}$ is supported on $[1/4,1)$ with  mass of $2^{-\alpha}$ at $1/4$.  The histogram of the ratio $\left(\lambda_{(2)}/\lambda_{(1)}\right)^{2/\alpha}$ based on 1000 replications from the model \eqref{ex:4.6} with noise given by a $t$-distribution with $\alpha=1.5$ degrees of freedom, $n=1000$ and $p=200$ is displayed in Figure \ref{fig:1}.
\begin{figure} [h]
\begin{center}
\includegraphics[scale=.35]{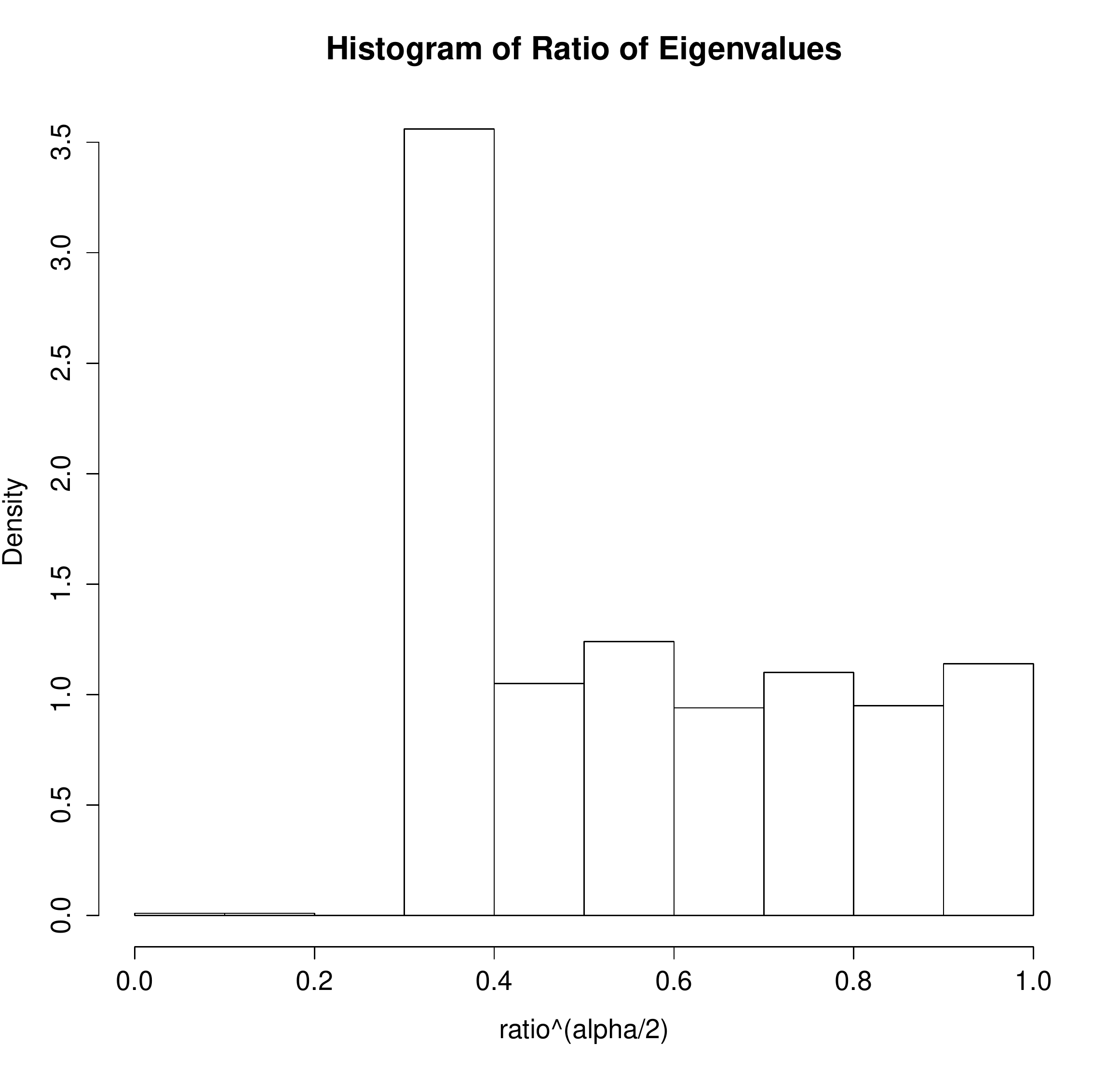}
\end{center}
\caption{Histogram based on 1000 replications of $\left(\lambda_{(2)}/\lambda_{(1)}\right)^{2/\alpha}$ from model \eqref{ex:4.6}.}
\label{fig:1}
\end{figure} Observing that $2^{-\alpha}=0.3536\ldots$,
the histogram is remarkably close to what one would expect from a sample from the truncated
uniform \ds , $2^{-\alpha}\, \1_{\{U<2^{-\alpha}\}}+U\, \1_{\{U\ge 2^{-\alpha} \}}$.
The mass of the limiting discrete component of the ratio can be much
larger if one conditions on $a_{np}^{-2}\lambda_{(1)}$ being large.
Specifically, for any $\epsilon\in (0,1/4)$ and $x>0$,
\begin{equation*}
\lim_{n\to\infty} \P(\epsilon<\lambda_{(2)}/\lambda_{(1)}\le 1/4|\lambda_{(1)}>a_{np}^2x)
=\P(\Gamma_1/\Gamma_2\le 2^{-\alpha}|\Gamma_1<(x/8)^{-\alpha/2})=G(x)\,.
\end{equation*}
The function $G$ approaches $1$ as $x\to\infty$ indicating the speed at which the
two largest eigenvalues get linearly related; see Figure \ref{fig:2} for a graph of $G$ in the case $\alpha=1.5$.
\begin{figure} [h]
\begin{center}
\includegraphics[scale=.35]{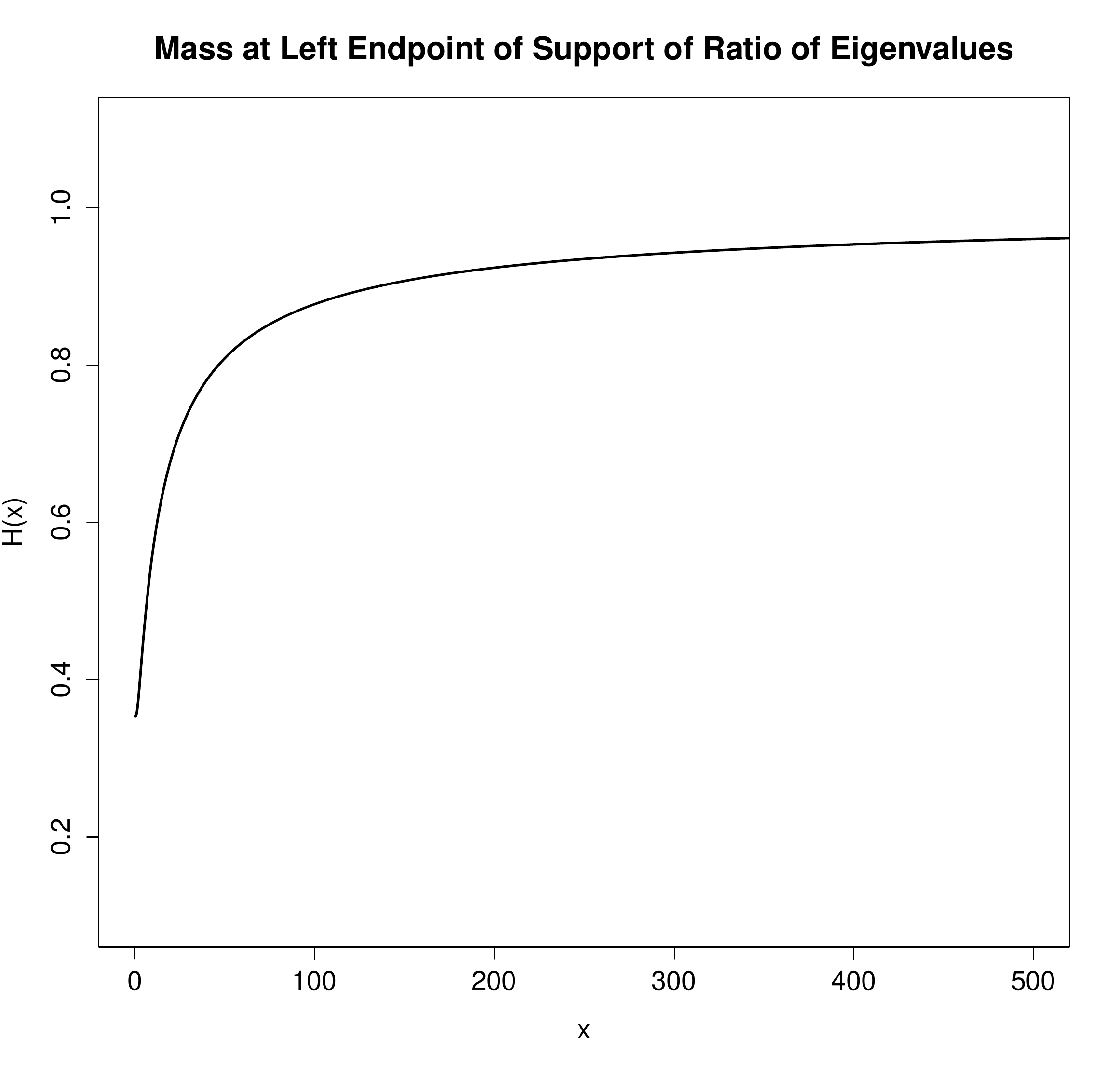}
\end{center}
\caption{Graph of $G(x)=\P(\Gamma_1/\Gamma_2\le2^{-\alpha}|\Gamma_1<(x/8)^{-\alpha/2})$ when $\alpha=1.5$.}
\label{fig:2}
\end{figure}
In addition, from Remark~\ref{rem:4.5}, we also have
\begin{equation*}
\dfrac{\la_{(1)}}{\la_1+\cdots +\la_p} \cid \dfrac{4}{5}\,
\dfrac{\Gamma_1^{-2/\alpha}}{\sum_{i=1}^\infty \Gamma_i^{-2/\alpha}}\,.
\end{equation*}
Clearly, the limit \rv\ is stochastically smaller than what one would get in the iid case; see \eqref{eq:limit}.
\end{example}
\begin{example}\label{exam:spectralgap}\rm
The previous example also illustrates the behavior of the two largest eigenvalues in the general case
when the rank $r$ of the matrix $\M$ is larger than one.
We have in general
\begin{equation*}
\dfrac{\lambda_{(2)}}{\lambda_{(1)}}\cid
\frac{v_2}{v_1}\,\1_{\{U<(v_2/v_1)^{\alpha/2}\}}+U^{2/\alpha}\,\1_{\{U\ge (v_2/v_1)^{\alpha/2}\}}\,.
\end{equation*}
In particular, the limiting
{\em self-normalized spectral gap} has \rep
\beao
\dfrac{\la_{(1)}-\la_{(2)}}{\la_{(1)}} \cid \dfrac{v_1-v_2}{v_1}\,\1_{\{U<(v_2/v_1)^{\alpha/2}\}}+(1-U^{2/\alpha})\,\1_{\{U\ge (v_2/v_1)^{\alpha/2}\}}\,.
\eeao
The limiting variable assumes values in $(0,1-v_2/v_1]$  and has an atom at the right end-point.
This is in contrast to the iid case
and to the case when $r=1$ (hence $v_2=0$) including the case of iid rows and the separable case; see Example~\ref{exam:separable}.
\end{example}
\begin{example}\label{exam:separable}
{\em We consider the separable case when  $h_{kl}=\theta_kc_l$, $k,l \in \Z$, where $(c_l)$, $(\theta_k)$ are
  real sequences such that the conditions on $(h_{kl})$ in
  Theorem~\ref{thm:mains} hold. In this case,
\begin{equation*}
\M= \sum_{l\in \Z} c_l^2\;\; (\theta_i\theta_j)_{i,j \in \Z}\,.
\end{equation*}
Note that $r=1$ with the only non-negative eigenvalue
\beao
v_1=\sum_{l\in \Z} c_l^2\;\;\sum_{k\in \Z}\theta_k^2\,.
\eeao
In this
case, the limiting point process in Theorem~\ref{cor:1} is a PRM on
$(0,\infty)$ with mean measure of $(y,\infty)$ given by $(v_1/
y)^{\alpha/2}$, $y>0$. The normalized eigenvalues have similar asymptotic behavior as
in the case of iid entries. For example, the log-spacings have the same limit as in the iid case for fixed $k$,
\beao
\big(\log \la_{(1)}-\log \la_{(2)},\ldots,\log \la_{(k+1)}-\log \la_{(k)}\big)\std
-\dfrac{2}{\alpha}\,\big(\log(\Gamma_1/\Gamma_2),\ldots,\log (\Gamma_{k}/\Gamma_{k+1})\big)\,.
\eeao
The same observation applies to the ratio of the largest eigenvalue and the trace in the case $\alpha\in (0,2)$:
\beao
\dfrac{\la_{(1)}}{{\rm tr} (\X \X')}=\dfrac{\la_{(1)}}{\la_1+\cdots +\la_p} \std
\dfrac{\Gamma_1^{-2/\alpha}}{\sum_{i=1}^\infty \Gamma_i^{-2/\alpha}}\,.
\eeao
We also mentioned in Example~\ref{exam:spectralgap} that the \ds al limit of the self-normalized spectral gap has no atom
as in the iid case.
}
\end{example}

\subsection{S\&P 500 data}\label{sec:sp500}
We conduct a short analysis of the largest eigenvalues of the univariate log-return time series which compose
the S\&P 500 stock index; see Section~\ref{subsec:1.2} for a description of the data.
Although there is strong empirical evidence that these univariate series have power-law tails
(see Figure~\ref{fig:SP500_tail_indices}) we  do not expect that they
have the same tail index. One way to proceed would be to ignore this fact because the tail indices are
in a close range and the differences are due to large sampling errors for estimating such quantities.  One could also collect \ts\ with similar
tail indices in the same group.
In this case, the dimension $p$ decreases. This grouping would be a rather arbitrary classification method.
We have chosen a third way: to use rank transforms.
This approach has its merits because it aims at standardizing the tails but it also has a major disadvantage: one destroys
the covariance structure underlying the data.
\par
Given a $p\times n$ matrix $(R_{it})_{i = 1, \cdots, p;t=1,\cdots, n}$,
we construct a matrix $\X$ via the rank transforms
\begin{equation*}
  X_{it} = -\Big[
    \log \Big(\frac{1}{n+1} \sum_{\tau=1}^n \1_{\{R_{i\tau} \leq R_{it}\}} \Big)
  \Big]^{-1} \,,\qquad i=1,\ldots,p;t=1,\ldots,n\,.
\end{equation*}
\begin{figure}[htb!]
  \centering
    \includegraphics[scale=0.550]{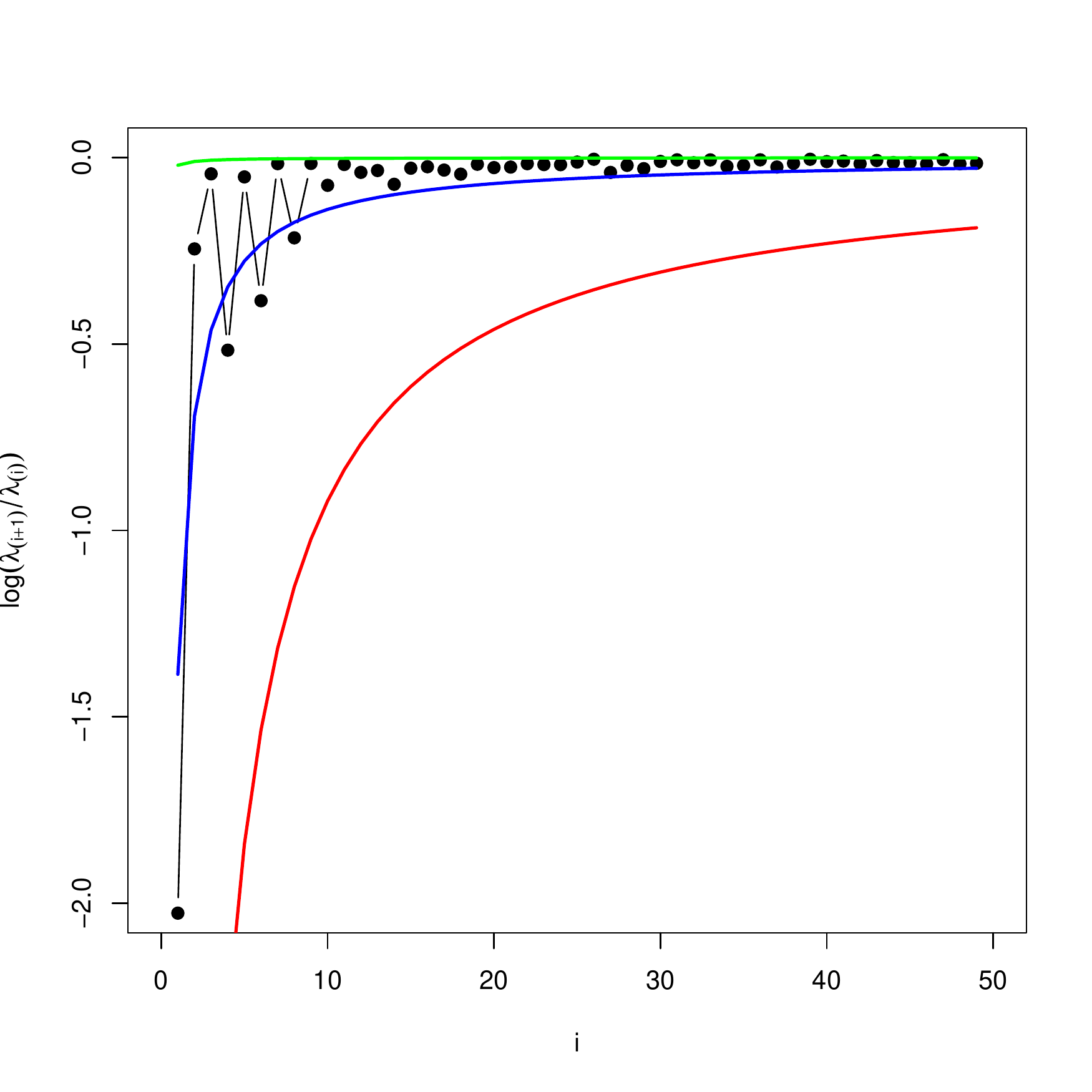}
  \caption{The logarithms of the ratios $\lambda_{(i+1)} / \lambda_{(i)}$ for the S\&P 500 series after rank transform.
We also show the 1, 50 and 99\% quantiles (bottom, middle, top lines, respectively) of the variables
$\log((\Gamma_i/ \Gamma_{i+1})^{2})$. }
  \label{fig:EigenRatio}
\end{figure}
\begin{figure}[htb!]
  \centering
    \includegraphics[scale=0.550]{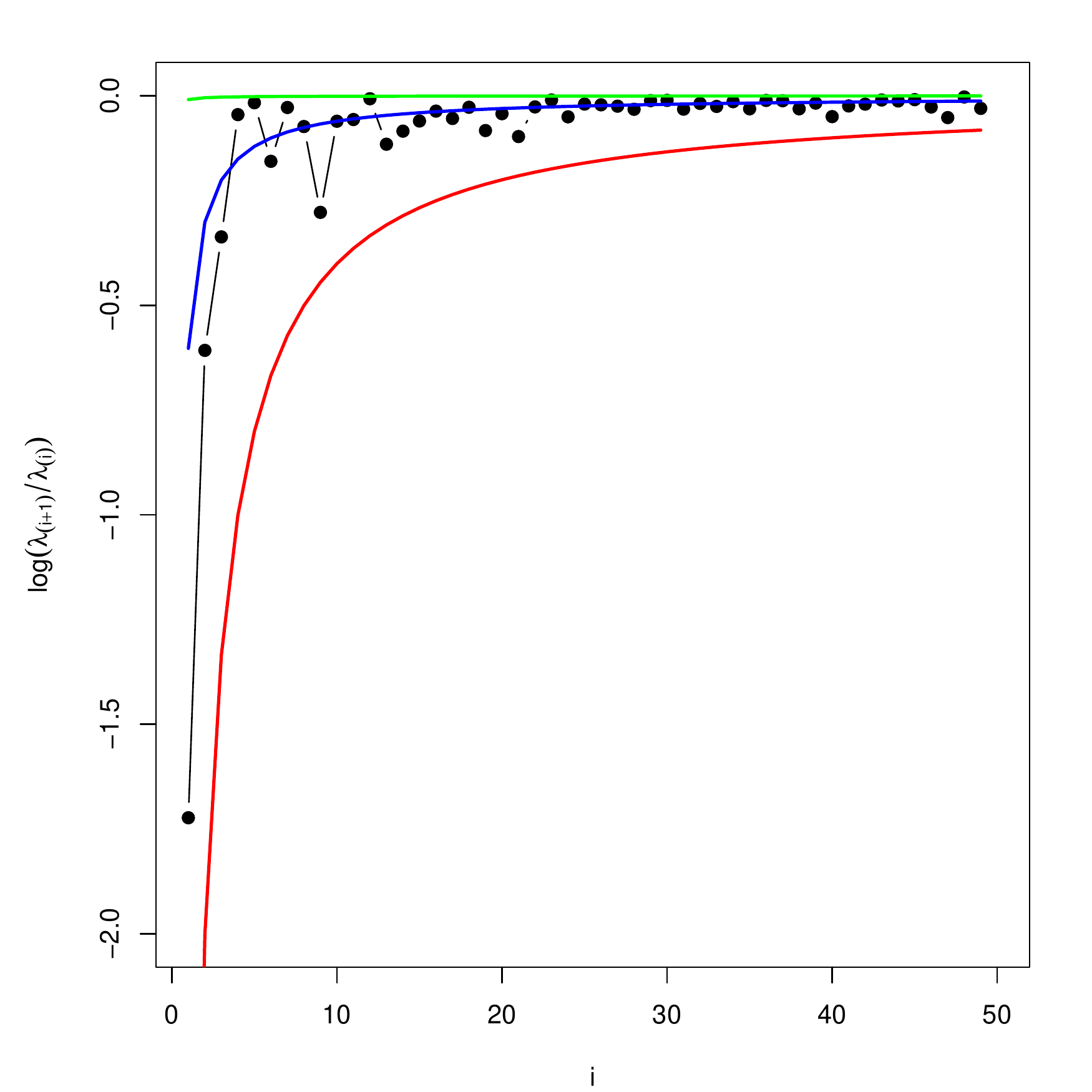}
  \caption{The logarithms of the ratios $\lambda_{(i+1)} / \lambda_{(i)}$ for the original (non-rank transformed) S\&P 500 log-return data.
We also show the 1, 50 and 99\% quantiles (bottom, middle, top lines, respectively) of the variables
$\log((\Gamma_i/ \Gamma_{i+1})^{2/2.3})$; see also Figure~\ref{fig:EigenRatio} for comparison.}\label{fig:22}
\end{figure}

%\begin{figure}[htb!]
%  \centering
%    \includegraphics[scale=0.550]{EigenRatioSP500_lin2_250_shown.pdf}
%  \caption{{\green Same as Figure \ref{fig:EigenRatio} but on original scale, no logarithm. Choose one.}}
%\end{figure}

If the rows $R_{i1}, \ldots, R_{in}$ were iid (or, more generally, stationary ergodic)
with a continuous \ds\  then
the averages under the logarithm would be  asymptotically uniform
on $(0,1)$ as $n \to \infty$. Hence $X_{it}$ would  be \asy ally standard
Fr\'echet $\Phi_1$-distributed.
In what follows, we assume that the aforementioned univariate \ts\ of the S\&P 500 index have undergone the rank transform and that
their marginal \ds s are close to $\Phi_1$; we
always use the symbol $\X$ for the resulting multivariate series.
\par
In Figure~\ref{fig:EigenRatio} we show the ratios of the consecutive ordered eigenvalues $(\la_{(i+1)}/\la_{(i)})$
of the matrix $\X\X'$. This graph shows the rather surprising fact that the ratios are close to one even for small values $i$.
We also show the 1, 50 and 99 \% quantiles of the variables $((\Gamma_{i}/\Gamma_{i+1})^{2/\alpha})$ calculated from the formula
\beam\label{eq:iid}
\P\big((\Gamma_{i} /\Gamma_{i+1})^{2/\alpha}\le x \big) = x^{i \cdot \alpha/2}, \quad x\in (0,1)\,.
\eeam
For increasing $i$, the \ds\ is concentrated closely to 1, in agreement with the \slln\ which yields
$\Gamma_{i} /\Gamma_{i+1}\stas 1$ as $i\to\infty$.
The \asy\ \ds s \eqref{eq:iid} correspond to the case when the matrix
$\M$ has rank $r=1$. It includes the iid and separable cases;
see Example~\ref{exam:separable}. The shown \asy\ quantiles are in agreement with the rank $r=1$ hypothesis.
\par
For comparison, in Figure~\ref{fig:22} we also show the ratios  $(\la_{(i+1)}/\la_{(i)})$ for the 
non-rank transformed S\&P 500 data and the 1, 50 and 99\% quantiles of the variables
$\log((\Gamma_i/ \Gamma_{i+1})^{2/\alpha})$, where we choose $\alpha=2.3$ motivated by the estimated tail indices in 
Figure~\ref{fig:SP500_tail_indices}.
The two graphs in Figure~\ref{fig:EigenRatio} and Figure~\ref{fig:22} are quite similar but the smallest ratios for the original data are slightly larger than
for the rank-transformed data.

%--------------------------------------------------------------------------------------
\subsection{Sums of squares of sample autocovariance matrices}\label{sec:possemidef}
In this section we consider some additive \fct s
of the squares of $\A_n(s)=\X_n(0)\X_n(s)'$ given by $\A_n(s)\A_n(s)'$ for $s=0,1,\ldots$. By definition of the singular values of a matrix
(see  \eqref{eq:sigma}), the non-negative definite
matrix $\A_n(s)\A_n(s)'$ has eigenvalues $(\la_i^2(s))_{i=1,\ldots,p}$.
\par
The following result is a corollary of Theorem~\ref{thm:mains}.
\begin{proposition}\label{thm:mainstr} Consider the linear process \eqref{eq:1} under
the conditions of Theorem~\ref{thm:mains}. Then the following statements hold for $s\ge 0$:
\begin{enumerate}
\item[$(1)$]
We consider two disjoint cases:
$\alpha \in (0,2)$ and $\beta\in (0,\infty)$, or
$\alpha\in [2,4)$ and $\beta$ satisfying \ref{Cbeta}. Then
\beao
a_{np}^{-4} \max_{i=1,\ldots,p} |\lambda_{(i)}^2(s)-\delta_{(i)}^2(s)| \cip 0, \quad \nto.
\eeao
\item[$(2)$]
Assume $\beta\in [0,1]$.
If $\alpha \in (0,2]$, $\E[Z^2]=\infty$ or $\alpha\in [2,4)$, $\E [Z^2]<\infty$ and $\beta \in (\alpha/2-1,1]$, then
\beao
a_{np}^{-4} \max_{i=1,\ldots,p} |\la_{(i)}^2(s)-(\gamma_{(i)}^\rightarrow(s))^2| \cip 0, \quad \nto.
\eeao
Assume $\beta>1$. If $\alpha \in (0,2]$, $\E[Z^2]=\infty$ or $\alpha\in [2,4)$, $\E [Z^2]<\infty$ and $\beta^{-1} \in (\alpha/2-1,1]$. Then
\beao
a_{np}^{-4} \max_{i=1,\ldots,p} |\la_{(i)}^2(s)-(\gamma_{(i)}^\downarrow(s))^2| \cip 0, \quad \nto.
\eeao
\end{enumerate}
\end{proposition}
To the best of our knowledge, sums of squares of sample autocovariance matrices were used first in the paper by Lam and Yao
\cite{lam:yao}; their \ts\ model is quite different from ours.
\begin{proof}
Part (1).
The proof follows from Theorem~\ref{thm:mains} if we can show that
\beao
a_{np}^{-2}\max_{i=1,\ldots,p} \big(\la_{(i)}(s)+\delta_{(i)}(s) \big)=O_\P(1)\,\quad \nto\,.
\eeao
We have by Theorem~\ref{cor:1},
\beam \label{eq:gdfg}
a_{np}^{-2}\max_{i=1\ldots,p} \la_{(i)}(s)=a_{np}^{-2}\la_{(1)}(s) \cid c\, \xi_{\alpha/2}\,,
\eeam
where $\xi_{\alpha/2}$ has a $\Phi_{\alpha/2}$ \ds . In view of Theorem~\ref{thm:mains}(1) we also have
\beao
a_{np}^{-2}\max_{i=1\ldots,p} \delta_{(i)}(s)\cid c\, \xi_{\alpha/2}\,.
\eeao
Therefore, again using Theorem~\ref{thm:mains}(1), we have
\beao
\lefteqn{a_{np}^{-4} \max_{i=1,\ldots,p} |\lambda_{(i)}^2(s)-\delta_{(i)}^2(s)|}\\
&\le &\big[a_{np}^{-2} \max_{i=1,\ldots,p} |\lambda_{(i)}(s)-\delta_{(i)}(s)|\big]\,
\big[a_{np}^{-2} \max_{i=1,\ldots,p}\big ( |\lambda_{(i)}(s)|+|\delta_{(i)}(s)|\big)\big]
\cip 0, \quad \nto.
\eeao
This proves part (1).\\[1mm]
Part (2). Now assume $\beta\in [0,1]$ and $\alpha \in (0,2]$, $\E[Z^2]=\infty$ or $\alpha\in [2,4)$, $\E [Z^2]<\infty$ and $\beta\in (\alpha/2-1,1]$. Then \eqref{eq:gdfg} is still true and we have by Theorem~\ref{thm:mains}(2) and Theorem~\ref{cor:1}
\beao
a_{np}^{-2}\max_{i=1\ldots,p} \gamma_{(i)}^{\rightarrow}(s)\cid c\, \xi_{\alpha/2}\,.
\eeao
We then have
\beao
\lefteqn{a_{np}^{-4} \max_{i=1,\ldots,p} |\lambda_{(i)}^2(s)-(\gamma_{(i)}^\rightarrow(s))^2|}\\
&\le &\big[a_{np}^{-2} \max_{i=1,\ldots,p} |\lambda_{(i)}(s) -\gamma_{(i)}^\rightarrow(s)|\big]\,
\big[a_{np}^{-2} \max_{i=1,\ldots,p}\big ( \lambda_{(i)}(s)+\gamma_{(i)}^\rightarrow(s)
\big)\big]
\cip 0\,, \qquad \nto.
\eeao
The proof of the remaining part is similar and therefore omitted.
\end{proof}
Now, using Proposition~\ref{thm:mainstr} and a continuous mapping argument, we can show
limit theory for the eigenvalues
\beao
w_{(1)}(s_0,s_1)\ge \cdots \ge w_{(p)}(s_0,s_1)\,,\qquad 0\le s_0\le s_1\,,
\eeao
of the non-negative definite random matrices
\begin{equation}\label{eq:sumA}
\sum_{s=s_0}^{s_1} \A_n(s) \A_n(s)'\,.
\end{equation}

\begin{proposition}\label{prop:sumsmal}
Assume $0\le s_0\le s_1$ and the conditions of Theorem~\ref{thm:mains} hold.
If $\alpha \in (0,4)$ and $\beta\in (0,1] \cap  (\alpha/2-1,1]$ then
\begin{equation*}
a_{np}^{-4} \max_{i=1,\ldots,p} |w_{(i)}(s_0,s_1)-\omega_{(i)}(s_0,s_1)| \cip 0, \quad \nto,
\end{equation*}
where $\omega_{(i)}(s_0,s_1)$ are the ordered values of the set $\{Z_{(i),np}^4 v_j(s_0,s_1), i=1,\ldots,p;j=1,2,\ldots\}$
and $(v_j(s_0,s_1))$ are the ordered eigenvalues of $ \sum_{s=s_0}^{s_1} \M(s)\M(s)'$.
\end{proposition}

\begin{example}\label{exam:additive}{\em
\begin{figure}[htb!]
  \centering
  \subfigure[]
{
    \includegraphics[scale=0.4]{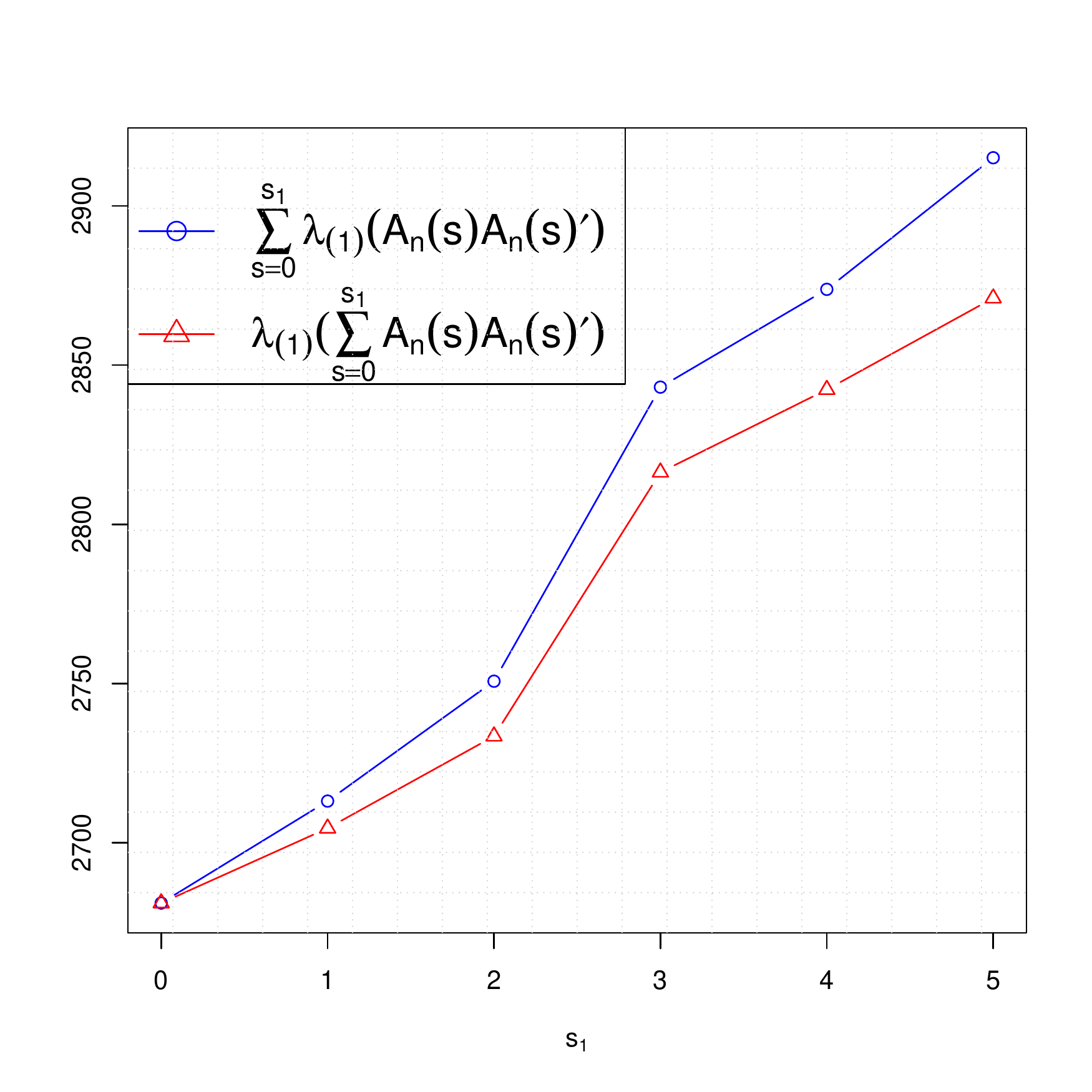}
    \label{fig:LamYao:a}
  }
  \subfigure[] {
    \includegraphics[scale=0.375]{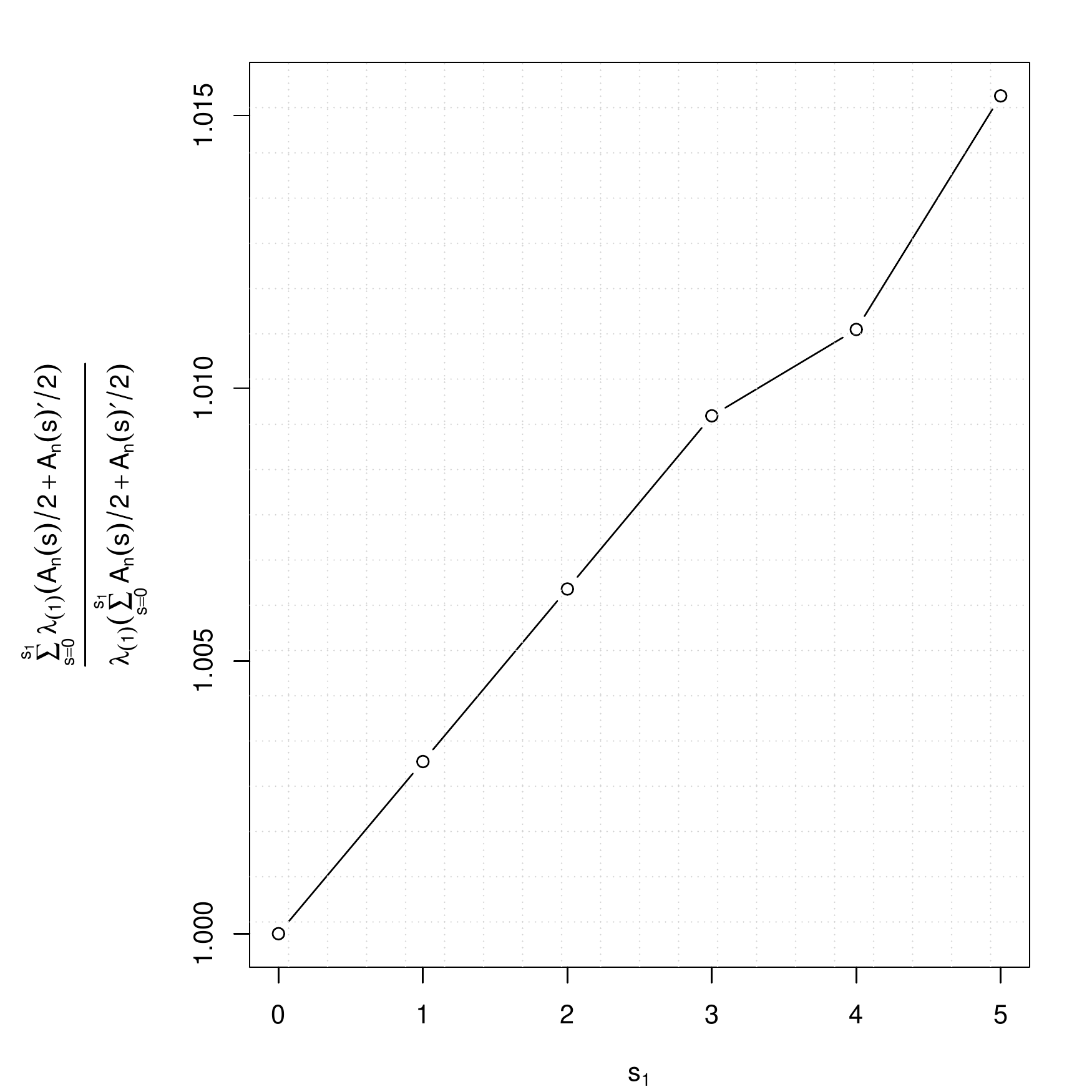}
    \label{fig:LamYao:b}
  }
  \caption{The largest eigenvalues of the sums of the squared autocovariance
    matrices compared with the sums of the largest eigenvalues of these matrices
    for the S\&P 500 data for different values $s_1$. The two values are surprisingly close to each other; mind the scale
of the $y$-axis. We also show their ratios.}
  \label{fig:LamYao}
\end{figure}
Recall the separable case from Example~\ref{exam:separable}, i.e.,
$h_{kl}=\theta_kc_l$, $k,l \ge 0$, where $(c_l)$, $(\theta_k)$ are
real sequences such that the conditions on $(h_{kl})$ in
Theorem~\ref{thm:mains} hold.
%Let $\alpha \in (0,2)$ and set
%\begin{equation*}
%\theta = (\theta_0, \theta_1,\theta_2,\ldots)' \quad \mbox{and} \quad
%c= (c_0,c_1,c_2, \ldots)'.
%\end{equation*}
Write $\Theta_{ij}=\theta_i \theta_j$. It is symmetric and has rank one; the only non-zero
eigenvalue is $\gamma_\theta(0)=\sum_{k=0}^\infty \theta_k^2$. Hence $\Theta$  is non-negative definite.
We get from \eqref{eq:m} that
\begin{equation*}
\M(s)=\gamma_c(s)\, \Theta, \quad s\ge 0\,,
\end{equation*}
where
\beao
\gamma_c(s)=\sum_{l=0}^\infty c_lc_{l+s}\,,\qquad s\ge 0\,.
\eeao
The matrix $\M(s)$ has the only non-zero eigenvalue $\gamma_c(s)\gamma_\theta(0)$.
The factors $(\gamma_c(s))$ can be positive or negative; they constitute the autocovariance \fct\ of a stationary linear process
with coefficients $(c_l)$.
Accordingly, $\M(s)$ is either non-negative or non-positive definite. %This explains why a symmetrization of $\M(s)$ is irrelevant.
Now we consider the non-negative definite matrix
\beao
\sum_{s=s_0}^{s_1} \M(s)\,\M(s)'= \sum_{s=s_0}^{s_1}\gamma_c^2(s)\,\Theta\Theta'\,.
\eeao
This matrix has rank $1$ and its largest eigenvalue is given by
\begin{equation*}
C_{c,\theta}(s_0,s_1)=\sum_{s=s_0}^{s_1}\gamma_c^2(s)\,\gamma_\theta^2(0)\,.
\end{equation*}
An application of
Proposition~\ref{prop:sumsmal} yields that the ordered eigenvalues of $a_{np}^{-4}\sum_{s=s_0}^{s_1}\A_n(s)\A_n(s)'$
are uniformly approximated by the quantities
\begin{equation}\label{eq:drtgdfg}
a_{np}^{-4} Z_{(i),np}^4 C_{c,\theta}(s_0,s_1)\,,\qquad  i=1,\ldots,p\,.
\end{equation}
Since
\beao
C_{c,\theta}(s_0,s_1)= \sum_{i=s_0}^{s_1} C_{c,\theta}(i,i)
\eeao
one gets the remarkable property that
\beao
&&a_{np}^{-4}\max_{i=1,\ldots,p} \Big| \la_{(i)}\Big(\sum_{s=s_0}^{s_1}\A_n(s)\A_n(s)'\Big)- Z_{(i),np}^4 C_{c,\theta}(s_0,s_1)\Big|\\
&=&a_{np}^{-4}\max_{i=1,\ldots,p} \Big|\sum_{s=s_0}^{s_1}\la_{(i)}(\A_n(s)\A_n(s)')- Z_{(i),np}^4 C_{c,\theta}(s_0,s_1)\Big|+o_P(1)\,.
\eeao
In particular, for $s_1\ge s_0$ we get the weak \con\ of the \pp es towards a PRM:
\beao
\sum_{i=1}^p \varepsilon_{a_{np}^{-4}\Big(\la_{i}\Big(\sum_{s=s_0}^{s_0}\A_n(s)\A_n(s)'\big)\,,\ldots,\la_{i}\big(\sum_{s=s_0}^{s_1}\A_n(s)\A_n(s)'\big)\Big)}
\std \sum_{i=1}^\infty \varepsilon_{\Gamma_i^{-4/\alpha} \Big(C_{c,\theta}(s_0,s_0),\ldots,C_{c,\theta}(s_0,s_1)\Big)}\,.
\eeao

}
\end{example}
\begin{example}\em
In Figure~\ref{fig:LamYao} we calculate the largest eigenvalues
$\la_{(1)}\big(\sum_{s=0}^{s_1}\A_n(s)\A_n(s)'\big)$ for $s_1=0,\ldots,5$ as well as the sums of the largest eigenvalues
$\sum_{s=0}^{s_1}\la_{(1)}\big(\A_n(s)\A_n(s)'\big)$
the log-return series from the S\&P 500 index described in Section~\ref{subsec:1.2}. The data are not rank-transformed.
We notice that
the two values are surprisingly close across the values $s_0=0,\ldots,5$. This phenomenon
could be explained by the structure of the eigenvalues in Example~\ref{exam:additive}.
Also note that the largest eigenvalue $\A_n(0)\A_n(0)'$ makes a major contribution to the values in Figure~\ref{fig:LamYao};
the contribution of the squares $\A_n(s)\A_n(s)'$, $s=1,\ldots,5,$ to the largest eigenvalue of the sum of squares is less substantial.
%Figure~\ref{fig:LamYao}(b) shows the same calculations for the matrices studied in Example~\ref{exam:4.6cont}. Here cancellation is possible (see \eqref{eq:Xrgsdf}) which is indicated by the fact that the red line does not follow the blue line as well as in (a).
\end{example}
%\section{\red Some conclusions?}\setcounter{equation}{0}
%{\em Say something about EVT and largest eigenvalues}

\appendix
\section{Auxiliary results}\label{appendix:A}\setcounter{equation}{0}

Let $(Z_i)$ be iid copies of $Z$ whose distribution satisfies
\begin{equation*}
\P(Z>x)\sim p_+ \dfrac{L(x)}{x^{\alpha}}\quad\mbox{and}\quad  \P(Z\le -x)\sim p_-
\dfrac{L(x)}{x^{\alpha}} \,,\quad \xto\,,
\end{equation*}
 for some tail index $\alpha>0$,
where $p_+,p_-\ge 0$ with $p_++p_-=1$ and $L$ is a slowly varying function. If $\E[|Z|]<\infty$ also assume $\E[Z]=0$. The product $Z_1Z_2$ is regular varying with the same index $\alpha$ and $\P(|Z_1Z_2|>x)= x^{-\alpha} L_1(x)$, where $L_1$ is slowly varying function different from $L$;
see Embrechts and Goldie \cite{embrechts:goldie:1980}.
Write
\begin{equation*}
S_n=Z_1+\cdots +Z_n\,,\quad n\ge 1,
\end{equation*} and consider a sequence $(a_n)$ such that $\P(|Z|>a_n)\sim n^{-1}$.

\subsection{Large deviation results}
The following theorem can be found in
Nagaev \cite{nagaev:1979} and Cline and Hsing
\cite{cline:hsing:1998} for $\alpha>2$ and $\alpha\le 2$,
respectively; see also  Denisov et al.~\cite{denisov:dieker:shneer:2008}.
\begin{theorem}\label{thm:nagaev}
Under the assumptions on the iid sequence $(Z_t)$
given above the following relation holds
\begin{equation*}
\sup_{x\ge c_n}\left|\dfrac{\P(S_n>x )}{n\P(|Z|>x)} -p_+ \right|\to 0\,,
\end{equation*}
where $(c_n)$ is any sequence satisfying $c_n/a_n\to  \infty$ for
$\alpha\le 2$ and $c_n\ge \sqrt{(\alpha-2)n\log n}$ for $\alpha>2$.
\end{theorem}

\subsection{A point process convergence result}
Assume that the conditions at the beginning of Appendix \ref{appendix:A} hold.
Consider a sequence of iid copies $(S_{n}^{(t)})_{t=1,2,\ldots}$
of $S_n$ and the sequence of point processes
\begin{equation*}
N_n= \sum_{t=1}^{p} \vep_{a_{np}^{-1} S_{n}^{(t)}}, \quad n=1,2,\ldots\,,
\end{equation*}
for an integer sequence $p=p_n\to\infty$. We assume that the state space of the
point processes $N_n$ is $\overline{\R}_0=[\R\cup\{\pm \infty\}]\backslash \{0\}$.
\begin{lemma}\label{lem:ppr}
Assume $\alpha \in (0,2)$ and the
conditions of  Appendix~\ref{appendix:A}
on the iid sequence $(Z_t)$ and the normalizing sequence $(a_n)$. Then the limit relation
$N_n\cid N$ holds in the space of point measures on $\overline{\R}_0$
equipped with the vague topology (see
\cite{resnick:1987,resnick:2007})
for a  Poisson random measure $N$ with state space $\overline{\R}_0$ and intensity measure $\mu_\alpha(dx)=\alpha |x|^{-\alpha-1} (p_+ \1_{\{x>0\}}+ p_- \1_{\{x<0\}}) dx$.
\end{lemma}
\begin{proof}
According to Resnick \cite{resnick:1987}, Proposition 3.21, we need to
show that
$p\, \P(a_{np}^{-1}S_n\in \cdot)\civ \mu_\alpha
$,
where $\civ$ denotes vague convergence of Radon measures on  $\overline{\R}_0$.
Observe that we have $a_{np}/a_n\to\infty$ as $\nto$. This fact and
$\alpha\in (0,2)$ allow one to apply
Theorem~\ref{thm:nagaev}:
\begin{equation*}
\dfrac{\P( S_n >x a_{np})}{n\,\P(|Z|>  a_{np})}\to p_+ x^{-\alpha}
\quad \mbox{and}\quad \dfrac{\P( S_n \le -x a_{np})}{n\,\P(|Z|>
  a_{np})}\to p_-\, x^{-\alpha}\,,\quad x>0\,.
\end{equation*}
On the other hand, $n\,\P(|Z|>  a_{np})\sim p^{-1}$ as $\nto$.
This proves the lemma.
\end{proof}

\section*{Acknowledgments}\setcounter{equation}{0}

We thank Olivier Wintenberger for reading the manuscript and fruitful discussions.

%\bibliography{libraryjohannes}

\end{document}